\documentclass[bj]{imsart}

\usepackage{amssymb}
\usepackage{amsfonts,amsmath,amsthm}
\usepackage{stackrel, subfigure}
\usepackage{graphicx}
\usepackage{pgf}

\usepackage{upgreek}

\usepackage{amssymb}
\usepackage{amsfonts}
\usepackage{mathrsfs}
\usepackage{upgreek}
\usepackage{graphicx}

\usepackage{appendix}

\providecommand{\U}[1]{\protect\rule{.1in}{.1in}}

\newtheorem{theorem}{Theorem}

\newtheorem{corollary}[theorem]{Corollary}

\newtheorem{lemma}[theorem]{Lemma}

\newtheorem{proposition}[theorem]{Proposition}

\newcommand{\ignore}[1]{}

\DeclareFontFamily{OT1}{pzc}{}
\DeclareFontShape{OT1}{pzc}{m}{it}{<-> s * [1.24] pzcmi7t}{}
\DeclareMathAlphabet{\mathpzc}{OT1}{pzc}{m}{it}

\begin{document}

%\tableofcontents

	\begin{frontmatter}
		
		\title{Posterior Probabilities: Nonmonotonicity, Asymptotic Rates, Log-concavity, and Tur\'{a}n's Inequality\protect\footnote{Previous version: July 2020 (Hebrew University of Jerusalem, Center for Rationality DP-736).}}
		\runtitle{Posterior Probabilities}

	\begin{aug}
		\author{\fnms{Sergiu}  \snm{Hart}\thanksref{a}\ead[label=e1]{hart@huji.ac.il}\ead[label=u1]{}},
		\and
		\author{\fnms{Yosef}  \snm{Rinott}\thanksref{b}%
			\ead[label=e3]{yosef.rinott@mail.huji.ac.il}%
			\ead[label=u1]{}}
		
		\runauthor{S. Hart and Y. Rinott}
		
		%\affiliation{Some University and Another University}
		
		\address[a]{The Hebrew  University   of Jerusalem, Federmann Center for the Study of Rationality, Einstein Institute of Mathematics, and Department of Economics. \printead{e1}}
		
		\address[b]{The Hebrew University of Jerusalem, Federmann Center for the Study of Rationality, and Department of Statistics. \printead{e3,u1}}
		
	\end{aug}

\begin{abstract}
	In the standard Bayesian framework data are assumed to be generated by a
	distribution parametrized by $\theta$ in a parameter space $\Theta$, over
	which a prior distribution $\pi$ is given. A Bayesian statistician quantifies
	the belief that the true parameter is $\theta_{0}$ in $\Theta$ by its
	posterior probability given the observed data. We investigate the behavior of
	the posterior belief in $\theta_{0}$ when the data are generated under some
	parameter $\theta_{1},$ which may or may not be the same as $\theta_{0}.$
	Starting from stochastic orders, specifically, likelihood ratio dominance,
	that obtain for resulting distributions of posteriors, we consider
	monotonicity properties of the posterior probabilities as a function of the
	sample size when data arrive sequentially. While the $\theta_{0}$-posterior is
	monotonically increasing (i.e., it is a submartingale) when the data are
	generated under that same $\theta_{0}$, it need not be monotonically
	decreasing in general, not even in terms of its overall expectation, when the
	data are generated under a different $\theta_{1}.$ In fact, it may keep going
	up and down many times, even in simple cases such as iid coin tosses. We obtain precise asymptotic rates when the data come
	from the wide class of exponential families of distributions; these rates
	imply in particular that the expectation of the $\theta_{0}$-posterior under
	$\theta_{1}\neq\theta_{0}$ is eventually strictly decreasing. Finally, we show
	that in a number of interesting cases this expectation is a log-concave
	function of the sample size, and thus unimodal. In the Bernoulli case we
	obtain this by developing an inequality that is related to Tur\'{a}n's
	inequality for Legendre polynomials.
	
\end{abstract}
		
		\begin{keyword}
			\kwd{Bayesian analysis}
			\kwd{stochastic and likelihood ratio orders}
			\kwd{sequential observations}
			\kwd{expected posteriors}
			\kwd{unimodality}
			\kwd{Legendre polynomials}	
			\kwd{exponential families}
		\end{keyword}
		
	\end{frontmatter}

\section{{Introduction\label{sec:int}}}

Consider a sequence of observations $x_{1},x_{2},...$ whose distribution is
governed by a parameter $\theta,$ in a standard Bayesian setup with a prior
distribution $\pi$ on the space of parameters $\Theta$. For simplicity assume
for now that $\Theta$ and the space of observations are finite (the results
extend readily to general observations and parameter spaces, as we will see later).
Let ${{\mathbb{P}}}_{\theta}$ denote the probability distribution under the
parameter $\theta,$ and let ${{\mathbb{P}}}=\sum_{\theta\in\Theta}\pi
(\theta){{\mathbb{P}}}_{\theta}$ denote the marginal probability; the
corresponding expectations are denoted by ${{\mathbb{E}}}_{\theta}$ and
${{\mathbb{E}}},$ respectively (thus ${{\mathbb{P}}}_{\theta}({{\mathbb{\cdot
			)\equiv P(}}}\cdot|\theta)$ and ${{\mathbb{{{\mathbb{E}}}}}}_{\theta
}\mathbb{[}{{\mathbb{\cdot]\equiv{{\mathbb{E}}}}}}[\cdot|\theta]).$

Let $\theta_{0}$ in $\Theta$ be a fixed value of the parameter. We are
interested in the way the belief in $\theta_{0}$ varies as one gets more and
more observations, i.e., as $n$ increases. We denote by $q_{n}^{\theta_{0}}$
the \emph{posterior} probability of $\theta_{0}$ at time $n;$ i.e., for every
sequence $s_{n}=(x_{1},...,x_{n})$ of observations up to time $n,$
\[
q_{n}^{\theta_{0}}\equiv q_{n}^{\theta_{0}}(s_{n})
{\;:=\;}
{{\mathbb{P}}}(\theta_{0}|s_{n})=\frac{{{\mathbb{P}}}_{\theta_{0}}(s_{n}%
	)\pi(\theta_0)}{{{\mathbb{P(}}}s_{n})}.
\]
As is well known, the sequence $q_{n}^{\theta_{0}}$ of posteriors is a
martingale with respect to the marginal probability ${{\mathbb{P}}}$, i.e.,
\begin{equation}
{{\mathbb{E[}}}q_{n+1}^{\theta_{0}}|s_{n}]=q_{n}^{\theta_{0}}(s_{n}%
)\label{eq:martingale}%
\end{equation}
for every $n$ and $s_{n}.$ Thus, given $s_{n},$ a new observation $x_{n+1}$
distributed according to ${{\mathbb{P}}}$ may increase or decrease the
posterior of $\theta_{0},$ but on average this posterior does not change.

Posterior probabilities are used in Bayesian hypothesis testing, where the
decision between hypotheses depends on their posterior probabilities (see
\cite{Berg}, Section 4.3.3). In sequential Bayes decision rules, monotonicity
properties of posteriors determine situations where more data lead to better
decisions vs. ones where more data could at times be misleading. This question
was raised and discussed in \cite{FK}. For a discussion of sequential Bayesian
inference and references see, e.g., \cite{Ferg} and \cite{Berg}. Setups where
more data may be harmful according to certain criteria of the
statistician are given, for example, in \cite{Barron} and \cite{Gru}. In
addition to the Bayesian setup, the questions addressed here are pertinent to
setups with agents that have different prior beliefs on the parameters. See
\cite{HR} for details, including applications to models with informed agents
vs. an uninformed market, and reputation-building models.

Suppose that the true parameter is $\theta_{0}$ (but this is, of course,
unknown to the observer). What can one say about the sequence of $\theta_{0}%
$-posteriors under ${{\mathbb{P}}}_{\theta_{0}}$? When the observations
$x_{n}$ are iid and the distributions ${{\mathbb{P}}}_{\theta}$ are distinct
(i.e., ${{\mathbb{P}}}_{\theta}\neq{{\mathbb{P}}}_{\theta^{\prime}}$ for
$\theta\neq\theta^{\prime},$ referred to as \textquotedblleft identifiable
parameters"), the Doob consistency theorem (see, e.g., \cite{VaV} and
\cite{Mill}) says that under ${{\mathbb{P}}}_{\theta_{0}}$ the sequence of
$\theta_{0}$-posteriors $q_{n}^{\theta_{0}}$ converges to $1$ almost surely
(i.e., except in a ${{\mathbb{P}}}_{\theta_{0}}$-null set). In fact, this
happens monotonically (see \cite{FK} and \cite{HR}, and the references
therein), in the sense that under ${{\mathbb{P}}}_{\theta_{0}}$ the posterior
of $\theta_{0}$ always increases on average; that is,
\begin{equation}
{{\mathbb{E}}}_{\theta_{0}}[q_{n+1}^{\theta_{0}}|s_{n}]\geq q_{n}^{\theta_{0}%
}(s_{n})\label{eq:ge-mart}%
\end{equation}
for every $n$ and $s_{n}$ (recall that ${{\mathbb{E}}}_{\theta_{0}}%
\mathbb{[}{{\mathbb{\cdot]\equiv{{\mathbb{E}}}}}}[\cdot|\theta_{0}]$ stands
for the expectation with respect to ${{\mathbb{P}}}_{\theta_{0}}$). This
submartingale inequality means that \emph{each additional observation
	increases on average the posterior of the true parameter, under the
	probability law of the true parameter.}

Taking the posterior probability of $\theta_{0}$ as one's belief in the model
determined by $\theta_{0}$, the expected belief in $\theta_{0}$ increases with
more data generated under $\theta_{0}$. Thus, as stated in \cite{FK}, Bayesian
inference does not lead one astray on average. In \cite{HR} we show that this
is in fact a consequence of an even stronger result, which holds for any
observation (and thus, in particular, for $x_{n+1}$ after $s_{n}$): the
distribution of the $\theta_{0}$-posterior under ${{\mathbb{P}}}_{\theta_{0}}$
dominates the distribution of that same $\theta_{0}$-posterior under the
marginal probability ${{\mathbb{P}}},$ where the domination is in the
likelihood ratio order, which is a strengthening of the usual stochastic order; see
Sections \ref{sec:sing} and \ref{sec:monot}. We thus get the submartingale
inequality (\ref{eq:ge-mart}), from which it follows that the overall
expectation ${{\mathbb{E}}}_{\theta_{0}}[q_{n}^{\theta_{0}}]\equiv
{{\mathbb{E}}}_{\theta_{0}}[{{\mathbb{P}}}({\theta_{0}}|s_{n})]$ of the
$\theta_{0}$-posterior $q_{n}^{\theta_{0}}$ under the probability
$\mathbb{P}_{\theta_{0}}$ is an increasing function of the number of
observations $n.$

Now suppose that the true parameter, $\theta_{1},$ is different from
$\theta_{0}.$ In the above case of iid observations and identifiable
parameters, the Doob consistency theorem now says that under ${{\mathbb{P}}%
}_{\theta_{1}}$ the sequence of $\theta_{0}$-posteriors $q_{n}^{\theta_{0}}$
converges to $0$ almost surely. When $\theta_{0}$ and $\theta_{1}$ are the
only possible parameter values (i.e., $\Theta=\{\theta_{0},\theta_{1}\}),$
by (\ref{eq:martingale}), (\ref{eq:ge-mart}), and ${{\mathbb{P}}}$ being the
average of ${{\mathbb{P}}}_{\theta_{0}}$ and ${{\mathbb{P}}}_{\theta_{1}},$ it
immediately follows that ${{\mathbb{E}}}_{\theta_{1}}[q_{n+1}^{\theta_{0}%
}|s_{n}]\leq q_{n}^{\theta_{0}}(s_{n});$ thus, each additional observation
decreases on average the posterior of a \textquotedblleft false" parameter,
under the probability law of the true parameter. However, this seemingly
natural property need \emph{not} hold when there are more than two possible
values of $\theta$ (see \cite{HR} for a simple example, and Section
\ref{sec:nonmonmult} below).

We next turn to consider monotonicity as a function of $n$ of $\psi(n)%
%TCIMACRO{\TeXButton{:=}{{\;:=\;}}}%
%BeginExpansion
{\;:=\;}%
%EndExpansion
{{\mathbb{E}}}_{\theta_{1}}[q_{n}^{\theta_{0}}]\equiv{{\mathbb{E}}}%
_{\theta_{1}}[{{\mathbb{P}}}({\theta_{0}}|s_{n})]$, the expectation over
$s_{n}$ of the $\theta_{0}$-posterior $q_{n}^{\theta_{0}}$ under the
probability ${{\mathbb{P}}}_{\theta_{1}},$ where $\theta_{1}\neq\theta_{0}.$
Consider for concreteness the simple setup of iid coin tosses. By Doob's
consistency result, $\psi(n)\rightarrow0$ as $n\rightarrow\infty$. This
convergence to zero need not however be monotonic. Indeed, when the true
parameter $\theta_{1}$ is \textquotedblleft close" to $\theta_{0}$ in a
suitable sense (which we will quantify in Proposition \ref{cor:uni}), it is
natural for data under $\theta_{1}$ to strengthen the belief in $\theta_{0}$
at first (i.e., for small $n),$ and so for $\psi(n)$ to increase for small $n$
(as is indeed the case when $\theta_{1}=\theta_{0}$; see (\ref{eq:ge-mart})).
Eventually, however, $\psi(n)$ must approach zero. Once $\psi(n)$ starts
decreasing in $n$, suggesting that evidence against $\theta_{0}$ is mounting,
can $\psi(n)$ increase again? While there may well be particular
realizations after which the $\theta_{0}$-posterior goes up (i.e.,
$q_{n+1}^{\theta_{0}}>q_{n}^{\theta_{0}}),$ and perhaps even particular data
$s_{n}$ after which the $\theta_{0}$-posterior is expected to go up (i.e.,
${{\mathbb{E}}}_{\theta_{1}}[q_{n+1}^{\theta_{0}}|s_{n}]>q_{n}^{\theta_{0}%
}(s_{n})),$ the \emph{overall expectation} is expected to be well behaved and
continue to go down (i.e., ${{\mathbb{E}}}_{\theta_{1}}[q_{n+1}^{\theta_{0}%
}]\leq{{\mathbb{E}}}_{\theta_{1}}[q_{n}^{\theta_{0}}]).$ Perhaps surprisingly,
that is \emph{not} the case: we provide simple examples where $\psi(n)$ is
\emph{not} unimodal in $n$ and may have multiple local maxima; that is, it can
go down and then up many times. Thus, after $\psi(n)$ starts decreasing and
the statistician may begin to doubt that $\theta_{0}$ is the true parameter,
the increase in $\psi(n)$ with more observations strengthens the
statistician's wrong belief that $\theta_{0}$ is the true parameter; the new
observations \emph{do} \textquotedblleft lead one astray on average." While
this is not a knife-edge phenomenon and may indeed happen, we show that
certain natural assumptions rule it out: the expected posterior is eventually
decreasing, and even log-concave and thus unimodal.

We now summarize the results on the behavior of $\psi(n)\equiv{{\mathbb{E}}%
}_{\theta_{1}}[q_{n}^{\theta_{0}}]$:

\begin{enumerate}
	\item When $\theta_{1}=\theta_{0},$ the sequence $\psi(n)$ is increasing in
	$n$ (as stated above, this is a consequence of the likelihood ratio dominance,
	which in turn implies the inequality (\ref{eq:ge-mart}); see Sections
	\ref{sec:sing} and \ref{sec:monot}).
	
	\item When $\theta_{1}\neq\theta_{0},$ the sequence $\psi(n)$ may increase in
	some range of values of $n$ and decrease in others, and may, for example,
	decrease first, then increase, and then decrease again (which, as we will show in
	Section \ref{sec:normnorm}, can happen already in the case of iid normal
	observations with a normal prior). Moreover, a large number of modes (i.e.,
	local maxima) can occur in the simple case of iid Bernoulli coin tosses (see
	Section \ref{sec:nonmonmult}).
	
	\item The sequence $\psi(n)$ is asymptotically equivalent to $C\sqrt{n}%
	\,w^{n}$ as $n\rightarrow\infty$ (for constants $C>0$ and $0<w\leq1)$ in the
	case of exponential families of distributions (discrete and continuous) with a
	continuous prior, where $w=1$ when $\theta_{1}=\theta_{0},$ and $w<1$ when
	$\theta_{1}\neq\theta_{0};$ in the latter case $\psi(n)$ is therefore strictly
	decreasing from some $n$ on (see Section \ref{sec:asympt}).
	
	\item The sequence $\psi(n)$ is log-concave in $n,$ and thus unimodal, in a
	number of scenarios: iid coin tosses with a uniform prior, iid normal
	observations with a normal prior (in a wide region of parameters), and iid
	exponential observations with an exponential prior (see Section
	\ref{sec:loggcon}).
\end{enumerate}

The paper is organized as follows. In Section \ref{sec:sing} we discuss
various order relations between posteriors under different distributions. In
Section \ref{sec:monot} we extend these results to sequences of observations
and see that $\psi(n)$ is increasing when $\theta_{1}=\theta_{0}.$ In Section
\ref{sec:nonmonmult} we exhibit situations in which $\psi(n)$ for $\theta
_{1}\neq\theta_{0}$ is not unimodal in $n$, and quantify, for Bernoulli
observations, the notion of $\theta_{0}$ and $\theta_{1}$ being close together such that
an observation under $\theta_{1}$ increases the belief in $\theta_{0}$ on
average. Section \ref{sec:asympt} provides the asymptotic analysis of
$\psi(n),$ and Section \ref{sec:loggcon} deals with cases where the sequence
$\psi(n)$ is unimodal and even log-concave. For Bernoulli observations this is
obtained by proving in Section \ref{sec:revtur} a reversal of the reverse
Tur\'{a}n inequality for orthogonal polynomials, which is of independent
interest. The Appendix contains discussions on possible extensions of the
results as well as some technical details.

\section{Preliminaries: A single signal}
\label{sec:sing} 
In this section we summarize results on various order
relations for posterior distributions that appeared with a somewhat different
emphasis in \cite{HR}. We consider a signal $s$, which can be a single
observation as well as multiple ones. For simplicity we now consider discrete
random variables $s$ and finite parameter spaces $\Theta$. The extension to
continuous random variables and prior distributions is straightforward; see
Appendix \ref{sec:general} and Section \ref{sec:asympt}.

We use the following standard notions and notation. The random variables $x$
and $y$ satisfy $y\geq_{{{\mathrm{st}}}}x$ (\emph{stochastic order}; also
known as \emph{first-order stochastic dominance}) or $y\geq_{{{\mathrm{icx}}}%
}x$ (\emph{increasing convex order}) if ${{\mathbb{E[}}}f(y)]\geq
{{\mathbb{E[}}}f(x)]$ for every increasing function or increasing convex
function $f$, respectively. By ``increasing" or ``convex" we do not mean
\textquotedblleft strictly" unless otherwise stated. The probability law (distribution) of a random
variable $x$ is denoted by ${\cal{L}}(x)$ (formally, ${\cal{L}%
}(x)=\mathbb{P}\circ x^{-1})$. Instead of $y\geq_{{{\mathrm{st}}}}x$ we may
write ${\cal{L}}(y)\geq_{{{\mathrm{st}}}}{\cal{L}}(x)$.

The \emph{likelihood ratio order}, denoted by $y\geq_{{{\mathrm{lr}}}}%
x$\thinspace or ${\cal{L}}(y)\geq_{{{\mathrm{lr}}}}{\cal{L}}(x)$,
is said to hold if ${{\mathbb{P}}}(y=t)/{{\mathbb{P}}}(x=t)$ is increasing in
$t$, or equivalently, ${{\mathbb{P}}}(y=t^{\prime}){{\mathbb{P}}}%
(x=t)\geq{{\mathbb{P}}}(y=t){{\mathbb{P}}}(x=t^{\prime})$ for all $t^{\prime
}>t$. This is stronger than the stochastic order: ${\cal{L}}%
(y)\geq_{{{\mathrm{lr}}}}{\cal{L}}(x)$ implies ${\cal{L}}%
(y)\geq_{{{\mathrm{st}}}}{\cal{L}}(x)$. Moreover, ${\cal{L}%
}(y)\geq_{{{\mathrm{lr}}}}{\cal{L}}(x)$ implies ${\cal{L}}(y|y\in
A)\geq_{{{\mathrm{lr}}}}{\cal{L}}(x|x\in A),$ and hence ${\cal{L}%
}(y|y\in A)\geq_{{{\mathrm{st}}}}{\cal{L}}(x|x\in A),$ for any
measurable subset $A$ of the real line. In fact, the latter condition of
stochastic dominance for every $A$ is equivalent to ${\cal{L}}%
(y)\geq_{{{\mathrm{lr}}}}{\cal{L}}(x)$; see \cite{Shak}.

In the standard Bayesian setup we have a prior $\pi$ with support $\Theta,$
and so $\pi(\theta)>0$ for every $\theta\in\Theta,$ and a random variable $s$
whose distribution depends on $\theta$. The conditional distribution
${{\mathbb{P}}}(s|\theta)$, namely, the distribution of $s$ given $\theta$,
is denoted by ${{\mathbb{P}}}_{\theta},$ and its probability law by
${\cal{L}}_{\theta}$. We also consider the marginal probability (also called the \textquotedblleft
prior predictive probability")
${{\mathbb{P(}}}s)%
%TCIMACRO{\TeXButton{:=}{{\;:=\;}}}%
%BeginExpansion
{\;:=\;}%
%EndExpansion
\sum_{\theta\in\Theta}{{\mathbb{P}}}_{\theta}(s)\pi(\theta)$, and denote its
law by ${\cal{L}}$. Expectations with respect to ${{\mathbb{P}}}$,
${{\mathbb{P}}}_{\theta},$ and ${{\mathbb{P}}}_{\Gamma}(s)%
%TCIMACRO{\TeXButton{:=}{{\;:=\;}}}%
%BeginExpansion
{\;:=\;}%
%EndExpansion
{{\mathbb{P}}}(s|\Gamma)=\sum_{\theta\in\Gamma}\pi(\theta){{\mathbb{P}}%
}_{\theta}(s)/\sum_{\theta\in\Gamma}\pi(\theta)$, for a set of parameters
$\Gamma\subset\Theta,$ are denoted by ${{\mathbb{E}}}$, ${{\mathbb{E}}%
}_{\theta}$, and ${{\mathbb{E}}}_{\Gamma}$, respectively. We use the notation
\[
q^{\theta}\equiv q^{\theta}(s)%
%TCIMACRO{\TeXButton{:=}{{\;:=\;}}}%
%BeginExpansion
{\;:=\;}%
%EndExpansion
{{\mathbb{P}}}(\theta|s)
\]
for the \emph{posterior} probability of $\theta$ (the \textquotedblleft%
$\theta$\emph{-posterior}" for short) given $s.$ We will compare random
variables like $q^{\theta}$ under different distributions, such as
${{\mathbb{P}}}$ and ${{\mathbb{P}}}_{\theta}$.

We now summarize several simple results with short proofs. They are given in
\cite{HR} with more details, interpretations, and related references.

\begin{proposition}
	\label{prop:p1p2}(i)\textbf{\ }Let $P_{1}$ and $P_{2}$ be two probability
	measures on a measure space ${\mathcal{S}}$ such that $P_{1}\ll P_{2}$ (i.e.,
	$P_{2}(s)=0$ implies $P_{1}(s)=0),$ and let $r(s)=P_{1}(s)/P_{2}(s)$ be the
	likelihood ratio.\footnote{When
		the ratio is $0/0$ define $r(s)$ arbitrarily; this will not matter since it
		occurs on a null event for both $P_{1}$ and $P_{2}.$} Then
	\[
	{\cal{L}}_{P_{1}}(r)\geq_{{{\mathrm{lr}}}}{\cal{L}}_{P_{2}}(r),
	\]
	where ${\cal{L}}_{P_{i}}$ denotes the probability law with respect to
	$P_{i},$ and for any increasing function $f$,
	\begin{equation}
	{\cal{L}}_{P_{1}}(f(r))\geq_{{{\mathrm{lr}}}}{\cal{L}}_{P_{2}%
	}(f(r)).\label{eq:f(r)}%
	\end{equation}
	(ii) In the Bayesian setup, for every $\theta$ the posterior $q^{\theta
	}(s)\equiv{{\mathbb{P}}}(\theta|s)$ of $\theta$ satisfies
	\begin{equation}
	{\cal{L}}_{\theta}(q^{\theta})\geq_{{{\mathrm{lr}}}}{\cal{L}%
	}(q^{\theta}).\label{eq:iii}%
	\end{equation}
	
\end{proposition}

\begin{proof}
	(i) For every value $t$ of $r$ let $B:=\{s:\frac{P_{1}(s)}{P_{2}(s)}=t\}$;
	then $P_{1}(r=t)=\sum_{s\in B}P_{1}(s)=t\sum_{s\in B}P_{2}(s)=tP_{2}(r=t)$. It
	follows that $\frac{P_{1}(r=t)}{P_{2}(r=t)}=t$, which is an increasing
	function of $t$, and so we have the result for $r$ by definition. The result
	for increasing $f$ then follows readily (see \cite{Shak}, Theorem 1.C.8).
	
	(ii) Setting $P_{1}={{\mathbb{P}}}_{\theta}$ and $P_{2}={{\mathbb{P}}}$ we
	have $q^{\theta}(s)=\pi(\theta)\frac{P_{1}(s)}{P_{2}(s)}=\pi(\theta)r(s)$, and
	the result follows from (i) since $\mathbb{P}_{\theta}\ll\mathbb{P}.$
\end{proof}

\noindent\textbf{Remark. }In (i) the ratio $P_{1}(r=t)/P_{2}(r=t)=t$ is a
\emph{strictly} increasing function of $t$---unless $r$ is constant, which
happens only when $P_{1}\equiv P_{2}$ (and then $r\equiv1)$---and the
domination is therefore \emph{strict}. In (ii) this happens except when the
signal $s$ is completely uniformative, i.e., when the posterior is identical
to the prior: $q^{\theta}(s)=\pi(\theta)$ for all $s.$ The strict domination
implies that the resulting inequalities, such as (\ref{eq:E-up}) below, are
strict (this is pointed out in \cite{HR}).

Since likelihood ratio order implies stochastic order, (\ref{eq:f(r)}%
)\textbf{\ }implies ${\cal{L}}_{P_{1}}(f(r))\geq_{{{\mathrm{st}}}%
}{\cal{L}}_{P_{2}}(f(r)),$ and so for any increasing function $f$ we
have
\[
\sum_{s}P_{1}(s)f\left(  \frac{P_{1}(s)}{P_{2}(s)}\right)  \geq\sum_{s}%
P_{2}(s)f\left(  \frac{P_{1}(s)}{P_{2}(s)}\right)  .
\]
The quantity on the right-hand side is known (for convex $f$) as
$f$\emph{-divergence}. Similarly, the likelihood ratio order relation
\eqref {eq:iii} implies stochastic order; that is, ${{\mathbb{E}}}_{\theta
}[f(q^{\theta}(s))]\geq{{\mathbb{E[}}}f(q^{\theta}(s))]$ for  increasing $f$.
Moreover, this holds also when
conditioning on a set of values of the posterior (such as being, say, more
than $1/2)$:%
\begin{equation}
{{{\mathbb{E}}}}_{\theta}[f(q^{\theta}(s))|q^{\theta}\in A]\geq{{{\mathbb{E[}%
}}}f(q^{\theta}(s))|q^{\theta}\in A]\label{eq:st-f-A}%
\end{equation}
for  increasing $f$ and $A\subseteq\lbrack0,1].$ Thus the posterior of $\theta$ when the data $s$ are generated according to
${{\mathbb{P}}}_{\theta}(s)$ is stochastically larger than the posterior when
the data are generated under ${{\mathbb{P}}}(s)$. Taking $f(x)=x$ we obtain
\begin{equation}
{{\mathbb{E}}}_{\theta}[q^{\theta}]\geq{{\mathbb{E[}}}q^{\theta}]=\pi
(\theta),\label{eq:E-up}%
\end{equation}
where the inequality is strict for any informative signal $s$ (see the above remark), and the equality is the \emph{martingale property} of posteriors under
${{\mathbb{P}}}$ (see (\ref{eq:martingale})). Replacing the single parameter
$\theta_{0}$ with a set $\Gamma\subset\Theta$ of parameter values, the result
of \eqref {eq:iii} readily implies that
\begin{equation}
{\cal{L}}_{\Gamma}\left(  {{\mathbb{P}}}(\Gamma|s)\right)
\geq_{{{\mathrm{lr}}}}{\cal{L}}\left(  {{\mathbb{P}}}(\Gamma|s)\right)
\geq_{{{\mathrm{lr}}}}{\cal{L}}_{\Gamma^{c}}\left(  {{\mathbb{P}}%
}(\Gamma|s)\right)  \label{eq:iii2G}%
\end{equation}
(for the second $\geq_{{{\mathrm{lr}}}}$ use ${{\mathbb{P}}}=\pi
(\Gamma){{\mathbb{P}}}_{\Gamma}+\pi(\Gamma^{c}){{\mathbb{P}}}_{\Gamma^{c}}),$
and thus
\begin{equation}
{{\mathbb{E}}}_{\Gamma}[{{\mathbb{P}}}(\Gamma|s)]\geq{{\mathbb{E[P}}}%
(\Gamma|s)]=\pi(\Gamma)\geq{{\mathbb{E}}}_{\Gamma^{c}}[{{\mathbb{P}}}%
(\Gamma|s)].\label{eq:submartin1G}%
\end{equation}
Result (\ref{eq:iii2G}) is given in \cite{HR}, and \eqref {eq:submartin1G} is
given in \cite{FK}; see these papers for further results, references, and history.

When there are only two parameter values, say $\Theta=\{\theta_{0},\theta
_{1}\},$ (\ref{eq:iii2G}) becomes ${\cal{L}}_{\theta_{0}}(q^{\theta_{0}%
})\geq_{{{\mathrm{lr}}}}{\cal{L}}(q^{\theta_{0}})\geq_{{{\mathrm{lr}}}%
}{\cal{L}}_{\theta_{1}}(q^{\theta_{0}}).$ However, the latter dominance
relation need \emph{not }hold when there are additional parameter values in
$\Theta$; see, for instance, the example at the end of Section 1 in
\cite{HR}.\footnote{Where $\Theta=\{\alpha,\beta,\gamma\}$ and the dominance
	is reversed: ${\cal{L}}_{\gamma}(q^{\alpha})>_{{{\mathrm{lr}}}%
	}{{\mathcal{L(}}}q^{\alpha})$ (in the notation of the present paper, \cite{HR} shows that ${\cal{L}}_{\gamma}(1-q^{\alpha})<_{{{\mathrm{lr}%
	}}}{{\mathcal{L(}}}1{{\mathcal{-}}}q^{\alpha})).$} A case where it does hold
is provided in the proposition below, which applies, for instance, to the
Bernoulli and normal distributions, and many other exponential families  discussed later; see, e.g., \cite{Karlin} or \cite{Jew}. The parameter space is now an
interval on the real line, $\theta_{0}$ and $\theta_{1}$ are the two interval
ends, and the family of distributions ${{\mathbb{P}}}_{\theta}$ satisfies the
\emph{monotone likelihood ratio property} (\emph{MLRP}), i.e., ${{\mathbb{P}}%
}_{\theta^{\prime}}(s)/{{\mathbb{P}}}_{\theta}(s)$ is increasing in
$s\in{{\mathbb{R}}}$ for $\theta^{\prime}>\theta.$ The order $\geq
_{{{\mathrm{lr}}}}$ below can of course be replaced by the weaker
$\geq_{{{\mathrm{st}}}}$\thinspace or by comparisons of expectations.

\begin{proposition}
	\label{prop:sto}Let $\Theta=[\theta_{0},\theta_{1}]\subset{{\mathbb{R}}}$ and
	assume that ${{\mathbb{P}}}_{\theta}$ is a monotone likelihood ratio (MLRP)
	family. Then
	\[
	{\cal{L}}_{\theta_{0}}(q^{\theta_{0}})\geq_{{{\mathrm{lr}}}%
	}{\cal{L}}(q^{\theta_{0}})\geq_{{{\mathrm{lr}}}}{\cal{L}}%
	_{\theta_{1}}(q^{\theta_{0}}).
	\]
	
\end{proposition}

\begin{proof}
	The only new part is ${\cal{L}}(q^{\theta_{0}})\geq_{{{\mathrm{lr}}}%
	}{\cal{L}}_{\theta_{1}}(q^{\theta_{0}})$. First, $\frac{{{\mathbb{P}}%
		}_{\theta_{1}}(s)}{{{\mathbb{P}}}(s)}$ is increasing in $s$ by the MLRP
	assumption, because the denominator is a mixture of ${{\mathbb{P}}}_{\theta}%
	$'s over $\theta\leq\theta_{1}$, and so ${\cal{L}}_{\theta_{1}}%
	(s)\geq_{{{\mathrm{lr}}}}{\cal{L}}(s).$ Similarly, $q^{\theta_{0}%
	}(s)=\pi(\theta_{0})\frac{{{\mathbb{P}}}_{\theta_{0}}(s)}{{{\mathbb{P}}}(s)}$
	is decreasing in $s,$ and the argument in the proof of Proposition
	\ref{prop:p1p2} (i), now applied to a decreasing rather than increasing
	function and thus reversing the order, implies ${\cal{L}}(q^{\theta_{0}%
	})\geq_{{{\mathrm{lr}}}}{\cal{L}}_{\theta_{1}}(q^{\theta_{0}})$.
\end{proof}
We conclude with a simple symmetry between the $\theta_{0}$-posterior under
$\theta_{1}$ and the $\theta_{1}$-posterior under $\theta_{0}.$

\begin{proposition}
	\label{prop:symm}Let $\theta_{0}$ and $\theta_{1}$ be in $\Theta$. Then
	\[
	\frac{1}{\pi(\theta_{0})}{{\mathbb{E}}}_{\theta_{1}}[q^{\theta_{0}}]=\frac
	{1}{\pi(\theta_{1})}{{\mathbb{E}}}_{\theta_{0}}[q^{\theta_{1}}]=\frac{1}%
	{\pi(\theta_{0})\pi(\theta_{1})}{{\mathbb{E[}}}q^{\theta_{0}}\cdot
	q^{\theta_{1}}].
	\]
	
\end{proposition}

\begin{proof}
	We have
	\begin{align*}
	&\frac{1}{\pi(\theta_{0})}{{\mathbb{E}}}_{\theta_{1}}[q^{\theta_{0}}]  
	=\frac{1}{\pi(\theta_{0})}\sum_{s}\frac{{{\mathbb{P}}}_{\theta_{0}}%
		(s)\pi(\theta_{0})}{{{\mathbb{P}}}(s)}{{\mathbb{P}}}_{\theta_{1}}(s)\\
	  &=\frac{1}{\pi(\theta_{0})\pi(\theta_{1})}\sum_{s}\frac{{{\mathbb{P}}%
		}_{\theta_{0}}(s)\pi(\theta_{0})}{{{\mathbb{P}}}(s)}\frac{{{\mathbb{P}}%
		}_{\theta_{1}}(s)\pi(\theta_{1})}{{{\mathbb{P}}}(s)}{{\mathbb{P}}}(s)
  =\frac{1}{\pi(\theta_{0})\pi(\theta_{1})}{{\mathbb{E[}}}q^{\theta_{0}}\cdot
	q^{\theta_{1}}].
	\end{align*}
	The last expression is symmetric in $\theta_{0}$ and $\theta_{1},$ and so it is
	equal to $\frac{1}{\pi(\theta_{1})}{{\mathbb{E}}}_{\theta_{0}}[q^{\theta_{1}%
	}]$ as well.
\end{proof}

The same symmetry applies of course when we consider sequences of observations
(see for instance Corollary \ref{cor:uninn} below).

\section{Increasing posterior of the true state\label{sec:monot}}

We now consider observations that arrive sequentially and apply the results of
Section \ref{sec:sing} to obtain monotonicity and order relations as a
function of the sample size $n$.

The data consist of a process of observations $x_{1},x_{2},\ldots$ whose
distribution is ${{\mathbb{P}}}_{\theta}$, where $\theta$ lies in the
parameter space $\Theta$. At this point we make no assumptions about the
distributions of the observation process and the dependence structure (over
$n).$ We assume the standard Bayesian framework given in Section
\ref{sec:int}. Given the vector of observations up to stage $n,$ which we
denote by $s_{n}=(x_{1},\ldots,x_{n})$, the \emph{posterior} of $\theta$ at
time $n$ is $q_{n}^{\theta}\equiv{{\mathbb{P}}}\left(  \theta|s_{n}\right)  $.
Viewing $q_{n}^{\theta}$ as an $s_{n}$-measurable random variable, we obtain
that the sequence $q_{n}^{\theta}$ is a martingale with respect to the
probability ${{\mathbb{P}}},$ i.e., ${{\mathbb{E}}}[q_{n+1}^{\theta}%
|s_{n}]=q_{n}^{\theta}.$

Proposition \ref{prop:p1p2} (ii) applied to $x_{n+1}|s_{n}$ yields

\begin{proposition}
	\label{cor:nlr} For every $\theta,$ the posterior $q_{n+1}^{\theta}$ of
	$\theta$ at time $n+1$ satisfies
	\[
	{\cal{L}}_{\theta}\left(  q_{n+1}^{\theta}|s_{n}\right)  \geq
	_{{{\mathrm{lr}}}}{\cal{L}}\left(  q_{n+1}^{\theta}|s_{n}\right)
	.\label{eq:seqq}
	\]
	
\end{proposition}

Thus, given $s_{n},$ the $\theta$-posterior $q_{n+1}^{\theta}\equiv
{{\mathbb{P}}}(\theta|x_{n+1},s_{n})$ under the probability ${{\mathbb{P}}%
}_{\theta}(\cdot|s_{n})$ likelihood-ratio dominates that same $\theta
$-posterior under the probability ${{\mathbb{P}}}(\cdot|s_{n})$. Since, again,
$\geq_{{{\mathrm{lr}}}}$ implies $\geq_{{{\mathrm{st}}}},$ by Proposition
\ref{cor:nlr} and the martingale property of $q_{n}^{\theta}\equiv
{{\mathbb{P}}}\left(  \theta|s_{n}\right)  $ we get
\begin{equation}
{{\mathbb{E}}}_{\theta}[q_{n+1}^{\theta}|s_{n}]\geq{{\mathbb{E}}}%
[q_{n+1}^{\theta}|s_{n}]=q_{n}^{\theta}\label{eq:subsub}%
\end{equation}
(as in (\ref{eq:st-f-A}), one may
also condition on $q_{n+1}^{\theta}$ lying in a certain
set).
Proposition \ref{cor:nlr} generalizes Proposition \ref{prop:p1p2} (ii): given
any past data $s_{n}$, the posterior belief in $\theta$ with an additional
observation $x_{n+1}$ distributed according to ${{\mathbb{{{\mathbb{P}}}}}%
}_{\theta}$ likelihood-ratio dominates the same posterior when the additional
observation is distributed according to ${{\mathbb{P}}}$; \eqref {eq:subsub}
is then the consequent expectation comparison.

Inequality \eqref {eq:subsub} means that under ${{\mathbb{P}}}_{\theta}$ the
process $q_{n}^{\theta}$ is a submartingale. Since every increasing convex
(integrable) function of a submartingale is a submartingale, for any
increasing convex $f$ we have 
\begin{equation}
{{\mathbb{E}}}_{\theta}[f(q_{n+1}^{\theta})|s_{n}]\geq f(q_{n}^{\theta}).
\label{eq:gencor3}%
\end{equation}
Taking expectations on the two sides of \eqref {eq:subsub} and
\eqref {eq:gencor3} with respect to $s_{n}$ distributed under $\theta$ we obtain

\begin{corollary}
	\label{cor:mon}The expectation ${{\mathbb{E}}}_{\theta}[q_{n}^{\theta}]$ of
	the posterior probability of $\theta$ with respect to ${{\mathbb{P}}}_{\theta
	}$ is increasing in $n$, and, more generally, so is ${{\mathbb{E}}}_{\theta
	}[f(q_{n}^{\theta})]$ for any convex increasing function $f$.
\end{corollary}

Thus, with more data generated under $\theta$, the $\theta$-posterior
increases in the increasing convex order, and in particular in expectation.
The convergence of the expected posterior to $1$ (by Doob's theorem) is thus
monotone. We have obtained this result starting from a strong ordering: the
likelihood ratio order. That $q_{n}^{\theta}$ is a submartingale under
${{\mathbb{P}}}_{\theta}$ is shown in \cite{MS} and \cite{FK} (see also the
references therein), where other relevant results are given.

\section{Nonmonotonic and multimodal expected posteriors}

\label{sec:nonmonmult}

Assume now that the observations are generated under $\theta_{1}$, which is
different from $\theta_{0}$; then the posterior probability of $\theta_{0}$
given $s_{n}$ converges to zero as $n\rightarrow\infty$ by Doob's theorem.
This convergence may not be monotone, and in fact, if $\theta_{1}$ and
$\theta_{0}$ are close together, the expected posterior $\psi(n)\equiv\psi_{\theta
	_{0},\theta_{1}}(n)={{\mathbb{E}}}_{\theta_{1}}[q_{n}^{\theta_{0}}]$ may
increase as a function of $n$ for small $n$ before it starts decreasing. The
question that we address in this section is whether it is possible for the
expected $\theta_{0}$-posterior, with data generated under $\theta_{1}$, to
increase again after it starts decreasing. If some observations distributed
under $\theta_{1}$ cause the expected $\theta_{0}$-posterior to decrease, the
average belief in $\theta_{0}$ decreases, as it should under $\theta_{1}$, but
as we will show below it is possible for further observations generated under the
same $\theta_{1}$ to cause $\psi(n)$ to increase before it eventually
decreases to zero. Such an increase leads to an erroneous upturn of the
Bayesian statistician's degree of belief in $\theta_{0}$. Furthermore,
${{\mathbb{\psi(}}}n)$ need not be unimodal, and may fluctuate many times.

We present two examples of this behavior in the simplest case of iid coin
tosses, i.e., ${{\mathrm{Bernoulli}}}(\theta)$ observations. In Figure 1 the
sequence $\psi(n)$ decreases, then increases, and then decreases again to $0$.
Figure 2 provides a further counterexample to the unimodality of the sequence
${{\mathbb{\psi(}}}n)$, showing that it may have many
modes, and thus may alternate several times between increasing and decreasing. In both examples
the priors concentrate on three points $\theta_{j}$ (for $j=0,1,2)$ with
probabilities $\pi(\theta_{j})=\alpha_{j}$, and so the distribution of the
sufficient statistic $u_{n}:=\sum_{i=1}^{n}x_{i}$ is a mixture of three
${{\mathrm{Binomial}}}(n,\theta_{i})$ distributions: ${{\mathbb{P}}}%
(u_{n}=k)=\sum_{j=0}^{2}\alpha_{j}{{\mathbb{P}}}_{\theta_{j}}(u_{n}=k)$.
By continuity, it is clear that such examples are robust to small changes in the parameters and their associated probabilities, and priors having a larger or even continuous support with a similar behavior can
be constructed. The function $\psi_{\theta_{0},\theta_{1}%
}(n)={{\mathbb{E}}}_{\theta_{1}}[q_{n}^{\theta_{0}}],$ which is given by the
formula
\begin{equation}
\psi_{\theta_{0},\theta_{1}}(n)={{\mathbb{E}}}_{\theta_{1}}[q_{n}^{\theta_{0}%
}]=\sum_{k=0}^{n}\frac{{\binom{n}{k}}^{2}\theta_{0}^{k}(1-\theta_{0}%
	)^{n-k}\theta_{1}^{k}(1-\theta_{1})^{n-k}\pi(\theta_{0})}{{{\mathbb{P}}}%
	(u_{n}=k)}\label{eq:gg}%
\end{equation}
(and thus equals $\frac{\pi(\theta_{0})}{\pi(\theta_{1})}{{\mathbb{E}}%
}_{\theta_{0}}[q_{n}^{\theta_{1}}],$ as in Proposition \ref{prop:symm}), is
depicted in Figures 1 and 2. In Figure 3 we provide an example of iid
observations with a normal prior, which will be analyzed in Section
\ref{sec:normnorm}.

For some intuitive explanations, take Figure 1 first, where we consider the
posterior of $\theta_{0}=0.5,$ the data are generated under $\theta_{1}=0.65,$
and there is another possible parameter, $\theta_{2}=0.85.$ Initially the
expected posterior belief in $\theta_{0}$ decreases, as the observations under
$\theta_{1}$ make both $\theta_{1}$ and $\theta_{2}$ seem more likely (on
average). After $11$ observations the expected belief in $\theta_{0}$ starts
to increase, as $\theta_{2},$ whose distance from $\theta_1$ (under which the data are generated) is greater than that of  $\theta_{0}$, begins to seem less likely.
Eventually, after a further $70$ observations, the expected $\theta_{0}%
$-posterior begins its final descent to $0.$

Turning to Figure 2, we see that, in addition to long-term fluctuations similar to those
of Figure 1 (see also Figure 3), there are many short-term up and down
fluctuations\footnote{We have generated other examples where the number of
	modes exceeds $8$ by far.} that are likely due to the discreteness of the data
(cf. Figure 3) and of the time steps. In Proposition \ref{cor:uni} and
Corollary \ref{cor:uninn} below we try to shed some light on these
fluctuations. Interestingly, the fluctuations of the expected posterior
beliefs are rather sensitive to the values of the parameters; for example,
changing $\theta_{1}$ from $0.5$ to $0.45$, or $\theta_{2}$ from $0.85$ to
$0.8,$ yields a strictly decreasing sequence $\psi(n),$ with no up and down fluctuations.

\begin{figure}[!tbp]
	\centering
	\begin{minipage}[b]{0.48\textwidth}
		\includegraphics[width=\textwidth,height=6.5cm]{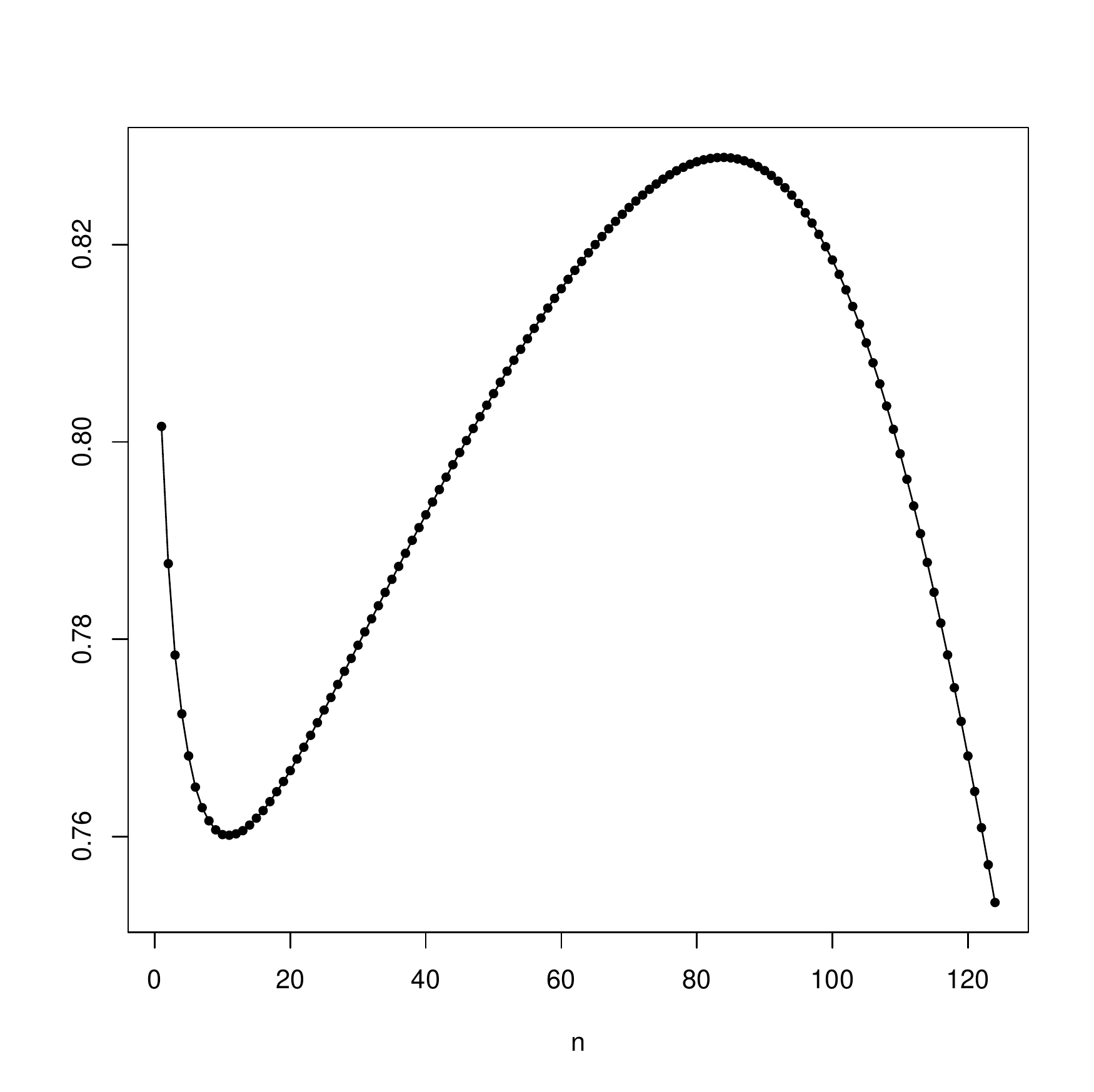}
		\caption{The sequence $\psi(n)=\mathbb{E}_{\theta_{1}}[q_{n}^{\theta_{0}}]$ for iid
			Bernoulli observations: $\Theta=\{\theta_{0},\theta_{1},\theta_{2}\}$ with
			$\theta_{j}=0.5,~0.65,~0.85,$ and prior probabilities $\alpha_{j}=\pi
			(\theta_{j})=4100/5001,~1/5001,~900/5001$ (for $j=0,1,2)$.}
	\end{minipage}
	\hfill 
	\begin{minipage}[b]{0.48\textwidth}
		\includegraphics[width=\textwidth,height=6.5cm]{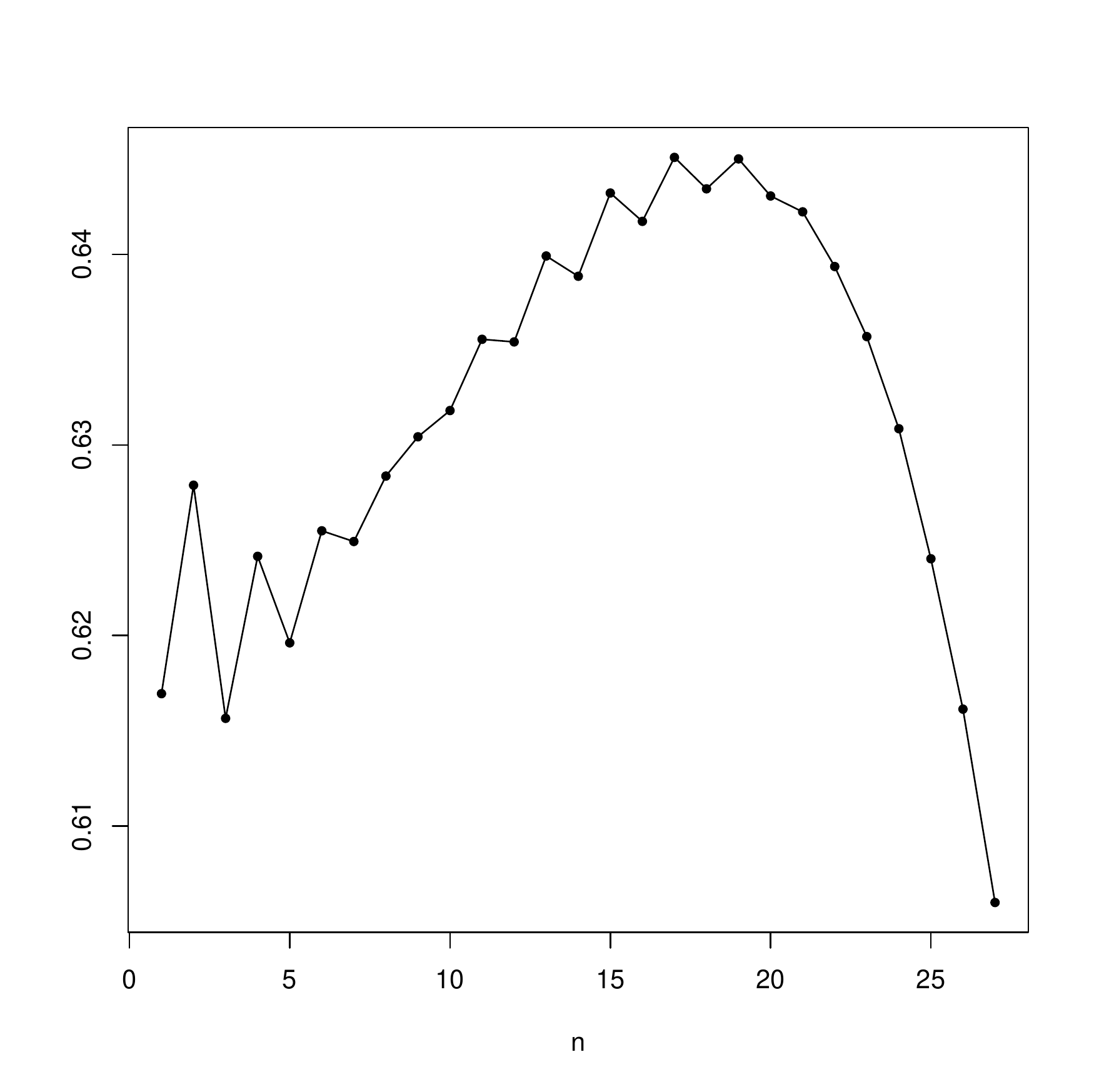}
		\caption{The sequence $\psi(n)=\mathbb{E}_{\theta_{1}}[q_{n}^{\theta_{0}}]$ for iid
			Bernoulli observations: $\Theta=\{\theta_{0},\theta_{1},\theta_{2}\}$ with
			$\theta_{j}=0.2,~0.5,~0.85,$ and prior probabilities $\alpha_{j}=\pi(\theta
			_{j})=2000/3001,~1/3001,~1000/3001$ (for $j=0,1,2);$ we see $8$ modes.}
	\end{minipage}
\end{figure}
\begin{figure}[h!]
	\centering
	\includegraphics[width=0.7\textwidth,height=6.5cm]{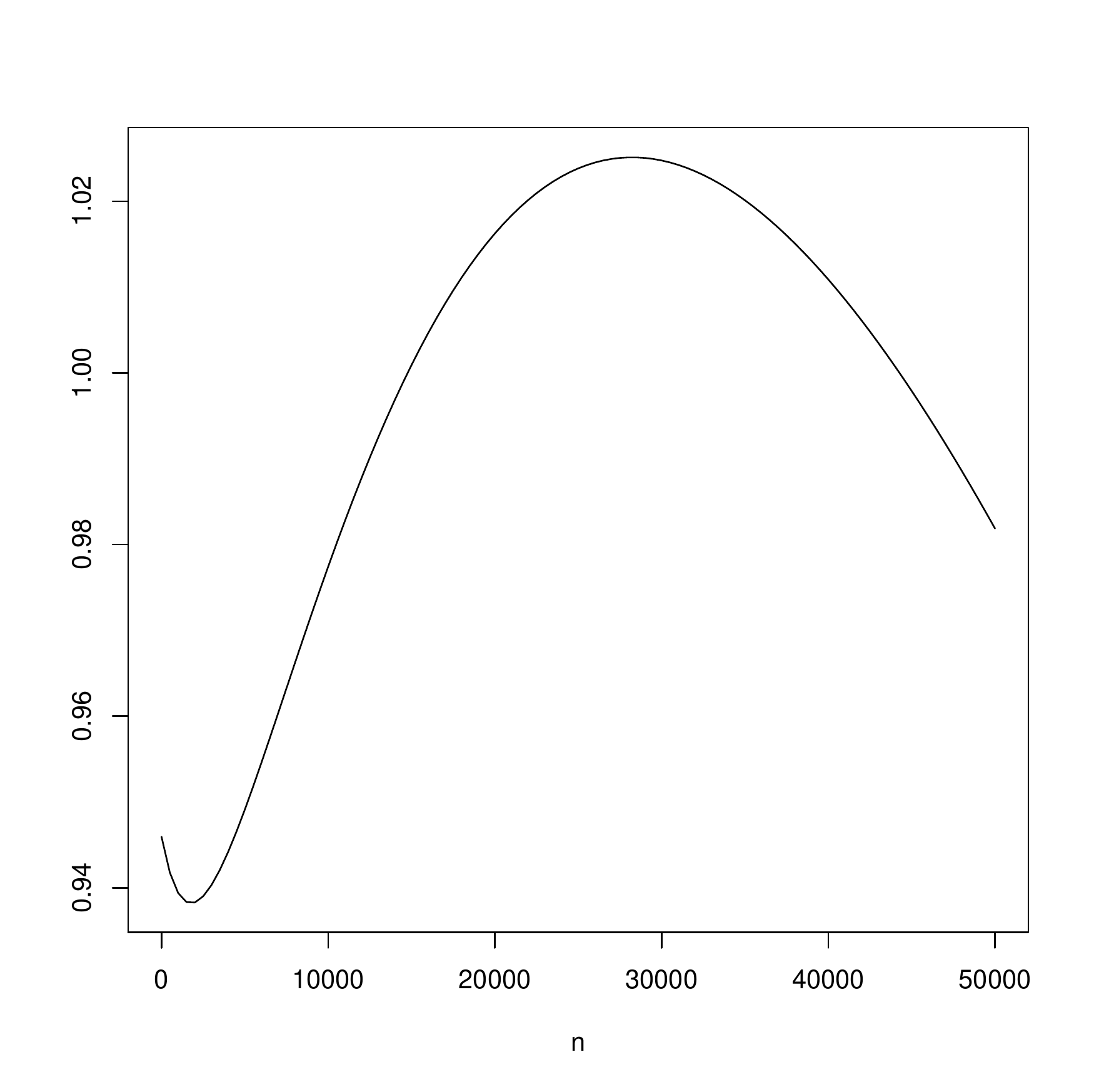}
	\caption{The sequence $\psi(n)=\mathbb{E}_{\theta_{1}}[q_{n}^{\theta_{0}}]$ for iid
		$\mathcal{N}(\theta,\sigma^{2})$ normal observations with a standard normal
		prior on $\theta$ (see (\ref{eq:nicenorm})): $\theta_{0}=-1/3,$~$\theta
		_{1}=1/3$ and $\sigma=100.$}
\end{figure}

In Figure 3 we have iid normal $\mathcal{N}(\theta,\sigma^{2})$ observations
with a large variance $(\sigma=100)$ and a standard normal prior (see Section
\ref{sec:normnorm}). For data generated under $\theta_{1}=1/3,$ the expected
belief in $\theta_{0}=-1/3$ starts by decreasing for nearly 2,000 periods,
following which it goes in the \textquotedblleft wrong direction" for a very
long time, increasing for about 30,000 additional observations; only then it
starts decreasing monotonically to its limit of $0.$ As we will show in Theorem
\ref{th:normal} (iv), additional up and down turns are not possible for normal distributions.

Consider now the effect of a single observation $x$ on the $\theta_{0}$-belief
when $x$ is generated under a different $\theta_{1}.$ As seen above, while in
many cases the belief in $\theta_{0}$ is reduced, this need not be so when
$\theta_{1}$ is close to $\theta_{0}$ and there are other points in $\Theta;$
indeed, the $\theta_{0}$-posterior strictly increases on average when
$\theta_{1}=\theta_{0}$ (see the remark after Proposition \ref{prop:p1p2}),
and thus, by continuity, also when $\theta_{1}$ and $\theta_{0}$ are close
enough together. We quantify precisely this \textquotedblleft closeness" of parameters
in the simplest case of a Bernoulli observation, as follows.

\begin{proposition}
	\label{cor:uni} Let $x$ be a ${{\mathrm{Bernoulli}}}(\theta)$ observation, and
	let $\pi$ be the prior probability of $\theta$ for $\theta\in\Theta
	\subseteq\lbrack0,1].$ Let $\overline{\theta}:={{\mathbb{E}}}[\theta
	]=\sum_{\theta\in\Theta}\theta\pi(\theta)$ be the (prior) average parameter.
	For any $\theta_{0}$ and $\theta_{1}$ in $\Theta,$ the inequality
	${{\mathbb{E}}}_{\theta_{1}}[q^{\theta_{0}}]\leq\pi(\theta_{0})$ holds if and
	only if $\overline{\theta}$ lies between $\theta_{0}$ and $\theta_{1},$ with
	equality if and only if $\overline{\theta}=\theta_{0}$ or $\overline{\theta
	}=\theta_{1}.$
\end{proposition}

The condition that $\overline{\theta}$ lies between $\theta_{0}$ and
$\theta_{1}$ is symmetric in $\theta_{0}$ and $\theta_{1},$ and so it is also
equivalent to ${{\mathbb{E}}}_{\theta_{0}}[q^{\theta_{1}}]\leq\pi(\theta_{1})$
(this equivalence follows from Proposition \ref{prop:symm} as well). Assume
therefore without loss of generality that $\theta_{0}\leq\theta_{1}.$
Proposition \ref{cor:uni} says that data under $\theta_{1}$ make the expected
posterior of $\theta_{0}$ lower than the prior of $\theta_{0}$ if and only if
$\theta_{0}\leq\overline{\theta}\leq\theta_{1},$ and make it higher than the
prior if and only if either $\theta_{0},\theta_{1}\leq\overline{\theta}$ or
$\theta_{0},\theta_{1}\geq\overline{\theta}$. Thus $\theta_{1}$ being
\textquotedblleft close" to $\theta_{0}$ in the sense that an observation
under $\theta_{1}$ increases the expected belief in $\theta_{0}$ (i.e.,
${{\mathbb{E}}}_{\theta_{1}}[q^{\theta_{0}}]\geq\pi(\theta_{0})$ holds) is
equivalent to $\theta_{1}$ being on the same side of $\overline{\theta}$ as
$\theta_{0}.$ In fact, this notion of \textquotedblleft belief-close" is
related to \textquotedblleft metric-close" as follows: the shorter the
distance $|\theta_{1}-\theta_{0}|$ between $\theta_{1}$ and $\theta_{0}$ is,
the fewer the priors $\pi$ yielding $\overline{\theta}$ in the interval
$[\theta_{0},\theta_{1}]$ there are, and so the fewer the priors $\pi$
yielding ${{\mathbb{E}}}_{\theta_{1}}[q^{\theta_{0}}]\geq\pi(\theta_{0})$
there are. Finally, when $\theta_{1}=\theta_{0},$ Proposition \ref{cor:uni}
gives (\ref{eq:E-up}), with strict inequality unless $\overline{\theta}%
=\theta_{0}.$

\begin{proof}
	Using ${{\mathbb{P}}}(x=1)=\overline{\theta}$ we have
	\begin{align}
	{{\mathbb{E}}}_{\theta_{1}}[q^{\theta_{0}}]  &  =\theta_{1}q^{\theta_{0}%
	}(1)+(1-\theta_{1})q^{\theta_{0}}(0)=\theta_{1}\frac{\theta_{0}\pi(\theta
		_{0})}{\overline{\theta}}+(1-\theta_{1})\frac{(1-\theta_{0})\pi(\theta_{0}%
		)}{1-\overline{\theta}}\nonumber\\
	&  =\pi(\theta_{0})V(\overline{\theta}), \label{eq:V}%
	\end{align}
	where
	\begin{equation}
	V(y){\;:=\;}\frac{\theta_{0}\theta_{1}}{y}+\frac{(1-\theta_{0})(1-\theta_{1}%
		)}{1-y}. \label{eq:pphhii}%
	\end{equation}
	The function $V$ is strictly convex and satisfies $V(\theta_{0})=V(\theta
	_{1})=1,$ and so $V(\overline{\theta})\leq1$ if and only if $\theta_{0}%
	\leq\overline{\theta}\leq\theta_{1}.$
\end{proof}

This extends to a sequence $x_{1},x_{2},\dots$ of iid ${{\mathrm{Bernoulli}}%
}(\theta)$ observations, where we obtain the condition for the expected
posterior to increase after the next observation, given the past data. Indeed,
using (\ref{eq:V}) for the observation $x_{n+1}$ that comes after data
$s_{n}=(x_{1},...,x_{n})$ yields

\begin{corollary}
	\label{cor:uninn} Let $x_{1},x_{2},\ldots$ be iid ${{\mathrm{Bernoulli}}%
	}(\theta)$ observations, let $\pi$ be the prior probability of $\theta$ for
	$\theta\in\Theta\subseteq\lbrack0,1],$ and denote by $\overline{\theta}%
	_{n}:={{\mathbb{E}}}[\theta|s_{n}]=\sum_{\theta\in\Theta}\theta{{\mathbb{P}}%
	}(\theta|s_{n})$ the average parameter conditional on $s_{n}.$ For any
	$\theta_{0}$ and $\theta_{1}$ in $\Theta$ we have
	\begin{equation}
	{{\mathbb{E}}}_{\theta_{1}}[q_{n+1}^{\theta_{0}}|s_{n}]=q_{n}^{\theta_{0}%
	}V(\overline{\theta}_{n})\;\;{{\text{and\ \ \ }}}{{\mathbb{E}}}_{\theta_{0}%
	}[q_{n+1}^{\theta_{1}}|s_{n}]=q_{n}^{\theta_{1}}V(\overline{\theta}%
	_{n}),\label{eq:V-n}%
	\end{equation}
	where the function $V$ is given by (\ref{eq:pphhii}), and so each one of the
	inequalities ${{\mathbb{E}}}_{\theta_{1}}[q_{n+1}^{\theta_{0}}|s_{n}]\leq
	q_{n}^{\theta_{0}}$ and ${{\mathbb{E}}}_{\theta_{0}}[q_{n+1}^{\theta_{1}%
	}|s_{n}]\leq q_{n}^{\theta_{1}}$ holds if and only if $\overline{\theta}_{n}$
	lies between $\theta_{0}$ and $\theta_{1}$.
\end{corollary}

Thus $\theta_{0}$ and $\theta_{1}$ are \textquotedblleft not close together" if they
are separated by $\overline{\theta}_{n}$, in which case an additional
observation under $\theta_{1}$ decreases the expected belief in $\theta_{0},$
and vice versa. The notion of closeness depends, of course, on $n$ and $s_{n}%
$. Since for data generated under $\theta_{1}$ (in the support of $\pi)$ we
have $\overline{\theta}_{n}\rightarrow\theta_{1}$, as $n$ increases
$\theta_{0}$ and $\theta_{1}$ may well fluctuate between being and not being close together given
$s_{n}$. Despite these fluctuations, we will see a large class of
natural setups (Theorem \ref{th:asymp} and Corollary \ref{c:asymp-monot})
where the overall expectation ${{\mathbb{E}}}_{\theta_{1}}[q_{n}^{\theta_{0}%
}]$ is strictly decreasing from some $n$ on.

\section{Asymptotic rates and eventual monotonicity}

\label{sec:asympt}

In this section we obtain the precise asymptotic behavior of ${{{\mathbb{E}}}%
}_{\theta_{1}}[q_{n}^{\theta_{0}}]$ as $n\rightarrow\infty$ in the commonly
used rich class of \emph{exponential families} of distributions; for clarity
we focus on one-dimensional families (see Appendix \ref{sec-a:asymp} for
extensions). The prior is now assumed to be continuous (which allows the use
of analysis tools), while the data may be discrete or continuous. All the
results of the previous sections are clearly seen to extend here, with
densities replacing probabilities as needed.

The general setup is as follows. Let ${{{\mathcal{X}}}}\subseteq{{{\mathbb{R}%
}}}$ be the space of observations, and let $\nu$ be a $\sigma$-finite measure on
the Borel sets of ${{{\mathcal{X}}}}$; for instance, take $\nu$ to be the
counting measure in the discrete case where ${{{\mathcal{X}}}}$ is a finite or
countable set, and the Lebesgue measure in the continuous case where
${{{\mathcal{X}}}}$ is a bounded or unbounded interval. Let $\Theta
\subseteq{{{\mathbb{R}}}}$ be the space of parameters, and let $\Pi,$ the
prior, be a probability measure on $\Theta;$ we assume that $\Theta$ is a
convex set and that $\Pi$ has a density $\pi(\theta)$ that is continuous and
strictly positive on $\Theta.$ Finally, the conditional-on-$\theta$
probability ${{\mathbb{P}}}_{\theta}(\cdot)\equiv{{\mathbb{P}}}(\cdot|\theta)$
on ${{\mathcal{X}}}$ has a density $p_{\theta}$ with respect to $\nu,$ given
by\footnote{Formally, $p_{\theta}$ is the Radon--Nikodym derivative
	${{\mathrm{d}}}{{{\mathbb{P}}}}_{\theta}/{{\mathrm{d}}}\nu$. In the discrete
	case where $\nu$ is the counting measure, $p_{\theta}(x)={{{\mathbb{P}}}%
	}_{\theta}(x)$ for all $x,$ and $\int_{Y}p_{\theta}(x)\,{{{\mathrm{d}}}}%
	\nu(x)=\sum_{x\in Y}p_{\theta}(x)=\sum_{x\in Y}{{\mathbb{P}}}_{\theta}(x)$ for
	every $Y\subseteq\mathcal{X}.$}
\begin{equation}
p_{\theta}(x)=\operatorname{exp}(\eta(\theta)T(x)-A(\eta(\theta
))-B(x))\label{eq:exp-fam}%
\end{equation}
for some functions $\eta,T,A,$ and $B,$ where $\eta$ is differentiable and
$\eta^{\prime}(\theta)>0$ for all $\theta\in\Theta.$ The following are well
known (see, e.g., \cite{Berk}): the function $A(\eta)$ is determined by the
functions $B$ and $T$ (use the condition $\int_{{{{\mathcal{X}}}}}p_{\theta
}(x)\,{{{\mathrm{d}}}}\nu(x)=1$ for every $\theta\in\Theta),$ it is infinitely
differentiable, $A^{\prime}(\eta(\theta))={{{\mathbb{E}}}}_{\theta}[T(x)],$
and $A^{\prime\prime}(\eta(\theta))={{{\mathbb{V}}}}ar_{\theta}(T(x)).$ We
assume that $A^{\prime\prime}(\eta(\theta))>0$ (that is, $T(x)$ is not
${{\mathbb{P}}}_{\theta}$-a.s. constant) for all $\theta\in\Theta.$ Let
$I(\theta):={{{\mathbb{E}}}}_{\theta}[(\partial/\partial\theta\log p_{\theta
}(x))^{2}]=-\mathbb{E}_{\theta}[\partial^{2}/\partial\theta^{2}\log p_{\theta
}(x)]$ be the \emph{Fisher information} at $\theta;$ for exponential families
(\ref{eq:exp-fam}) we have $I(\theta)=A^{\prime\prime}(\eta(\theta))\cdot
(\eta^{\prime}(\theta))^{2}={{{\mathbb{V}}}}ar_{\theta}(T(x))\cdot
(\eta^{\prime}(\theta))^{2}>0$.

Consider a sequence of iid observations $x_{1},x_{2},...$ and set
$s_{n}=(x_{1},...,x_{n})$ $\in{{{\mathcal{X}}}}^{n}.$ The density of the
$\theta_{0}$-posterior (for $\theta_{0}\in\Theta)$ is\footnote{We use the
	notation $p_{\theta}(\cdot)$ for the conditional-on-$\theta$ density of any
	variable; thus, $p_{\theta}(s_{n})=\prod_{i=1}^{n}p_{\theta}(x_{i})$.}
\[
q_{n}^{\theta_{0}}\equiv q_{n}^{\theta_{0}}(s_{n})=\frac{p_{\theta_{0}}%
	(s_{n})\pi(\theta_{0})}{p(s_{n})},
\]
where $p$ denotes the marginal density, i.e., $p(s_{n})=\int_{\Theta}%
p_{\theta}(s_{n})\pi(\theta)\,{{{\mathrm{d}}}}\theta.$ The $\theta_{1}%
$-expectation (for $\theta_{1}\in\Theta)$ of the $\theta_{0}$-posterior
is\footnote{The measure $\nu^{n}$ on ${{{\mathcal{X}}}}^{n}$ is the $n$-fold
	product of the measure $\nu$ on ${{{\mathcal{X}}}}.$}
\begin{equation}
\psi_{\theta_{0},\theta_{1}}(n)\equiv{{{\mathbb{E}}}}_{\theta_{1}}%
[q_{n}^{\theta_{0}}]=\int_{{{{\mathcal{X}}}}^{n}}\frac{p_{\theta_{0}}%
	(s_{n})\pi(\theta_{0})}{p(s_{n})}p_{\theta_{1}}(s_{n})\,{{{\mathrm{d}}}}%
\nu^{n}(s_{n}).\label{eq:psi}%
\end{equation}

In exponential families there is a simple relation between $\psi_{\theta
	_{0},\theta_{1}}$ and $\psi_{\theta,\theta}$ for an appropriate $\theta
\in\Theta.$

\begin{proposition}
	\label{p:theta-bar}Let $(p_{\theta})_{\theta\in\Theta}$ be a family of
	densities (\ref{eq:exp-fam}), and let $\pi$ be a positive density on the
	convex set $\Theta\subseteq\mathbb{R}.$ Let $\theta_{0},\theta_{1}\in\Theta.$
	Then
	\[
	\psi_{\theta_{0},\theta_{1}}(n)=\frac{\pi(\theta_{0})}{\pi(\theta_{2})}%
	\psi_{\theta_{2},\theta_{2}}(n)\,w^{n}%
	\]
	for every $n\geq1,$ where
	\begin{equation}
	\theta_{2}=\eta^{-1}\left(  \frac{\eta(\theta_{0})+\eta(\theta_{1})}%
	{2}\right)  \in\Theta\label{eq:theta2}%
	\end{equation}
	and\footnote{The constant $\sqrt{w}$ is known as the Bhattacharrya
		coefficient, and also as the Chernoff $1/2$-coefficient (see \cite{Bhat}, \cite{Cher}, or \cite{Nielsen} for a convenient reference).}
	\begin{equation}
	w=\left(  \int_{{{{\mathcal{X}}}}}\sqrt{p_{\theta_{0}}(x)\,p_{\theta_{1}}%
		(x)}\,{{{\mathrm{d}}}}\nu(x)\right)  ^{2}\leq1,\label{eq:w}%
	\end{equation}
	with equality (i.e., $w=1)$ if and only if $\theta_{0}=\theta_{1}.$
\end{proposition}

\begin{proof}
	The function $\eta$ is strictly increasing and continuous, and so it attains
	the value $(\eta(\theta_{0})+\eta(\theta_{1}))/2$ at a (unique) point between
	$\theta_{0}$ and $\theta_{1}$, and thus in $\Theta$ (which is a convex set);
	this is the point $\theta_{2}$ given by (\ref{eq:theta2}). From
	(\ref{eq:exp-fam}) we have
	\begin{equation}
	p_{\theta_{0}}(x)\,p_{\theta_{1}}(x)=w\,p_{\theta_{2}}(x)^{2},\label{eq:p0*p1}%
	\end{equation}
	for every $x$, where $w:=\operatorname{exp}(2A(\eta_{2})-A(\eta_{0}%
	)-A(\eta_{1}))$ and $\eta_{i}:=\eta(\theta_{i})$ (and so $\eta_{2}=(\eta
	_{0}+\eta_{1})/2).$ Therefore $p_{\theta_{0}}(s_{n})\,p_{\theta_{1}}%
	(s_{n})=w^{n}p_{\theta_{2}}(s_{n})^{2},$ yielding the result by (\ref{eq:psi}%
	). Formula (\ref{eq:w}) follows from (\ref{eq:p0*p1}): $\int_{{{{\mathcal{X}}%
	}}}\sqrt{p_{\theta_{0}}(x)\,p_{\theta_{1}}(x)}=\int_{{{{\mathcal{X}}}}}%
	\sqrt{w}\,p_{\theta_{2}}(x)=\sqrt{w}.$ Since $\int_{{{{\mathcal{X}}}}}%
	\sqrt{p_{\theta_{0}}(x)\,p_{\theta_{1}}(x)}\leq\int_{{{{\mathcal{X}}}}%
	}(p_{\theta_{0}}(x)+p_{\theta_{1}}(x))/2=1$ we get $w\leq1,$ with equality if
	and only if $p_{\theta_{0}}(x)=p_{\theta_{1}}(x)$ for all $x\in{{{\mathcal{X}%
	}}},$ which occurs if and only if $\theta_{0}=\theta_{1}$ (because $\eta$ is
	one-to-one and $T\ $is not constant). 
\end{proof}

Proposition \ref{p:theta-bar} thus reduces the analysis of the $\theta_{1}%
\neq\theta_{0}$ case to that of the $\theta_{1}=\theta_{0}$ case. The relation
(\ref{eq:p0*p1}), or, equivalently, $p_{\theta_{2}}=c\sqrt{p_{\theta_{0}%
	}\,p_{\theta_{1}}}$ (for the constant $c=1/\sqrt{w}$), says that the density
$p_{\theta_{2}}$ is proportional to the geometric average of the densities
$p_{\theta_{1}}$ and $p_{\theta_{2}};$ it is their \textquotedblleft
normalized geometric average." The exponential family (\ref{eq:exp-fam}) is
closed under this averaging operation by the convexity of $\Theta$. For
further discussions and extensions, see Appendix \ref{sec-a:asymp}.

When $\theta_{1}=\theta_{0}$ the sequence ${{{\mathbb{E}}}}_{\theta_{0}}%
[q_{n}^{\theta_{0}}]$ is monotonically increasing (by Corollary \ref{cor:mon}%
), whereas when $\theta_{1}\neq\theta_{0}$ the sequence ${{{\mathbb{E}}}%
}_{\theta_{1}}[q_{n}^{\theta_{0}}]$ converges to $0$ (by Doob's theorem).
Theorem \ref{th:asymp} strengthens these by yielding the precise asymptotic
rates. 

The standard notation \textquotedblleft$f(n)\sim g(n)$ as $n\rightarrow\infty
$\textquotedblright\ means \textquotedblleft asymptotic equivalence," i.e.,
$\lim_{n\rightarrow\infty}f(n)/g(n)=1$.

\begin{theorem}
	\label{th:asymp} Let $x_{1},x_{2},\ldots$ be iid observations distributed
	according to an exponential family with density $p_{\theta}$ given by
	(\ref{eq:exp-fam}), and let $\theta$ be distributed according to a prior
	having a continuous strictly positive density $\pi$ on a  convex set
	$\Theta\subseteq{{\mathbb{R}}}.$
	
	(i) Let $\theta_{0}$ be an interior point of $\Theta$; then
	\begin{equation}
	{{{\mathbb{E}}}}_{\theta_{0}}[q_{n}^{\theta_{0}}]\sim\frac{\sqrt{I(\theta
			_{0})}}{2\sqrt{\uppi }}\sqrt{n}{{\text{\ \ as }}}n\rightarrow\infty.
	\label{eq:asymp-theta0}%
	\end{equation}

	(ii) Let $\theta_{1}\neq\theta_{0}$ be in $\Theta$; then
	\begin{equation}
	{{{\mathbb{E}}}}_{\theta_{1}}[q_{n}^{\theta_{0}}]=\frac{\pi(\theta_{0})}%
	{\pi(\theta_{2})}{{{\mathbb{E}}}}_{\theta_{2}}[q_{n}^{\theta_{2}}]\,w^{n}%
	\sim\frac{\pi(\theta_{0})}{\pi(\theta_{2})}\frac{\sqrt{I(\theta_{2})}}%
	{2\sqrt{\uppi}}\sqrt{n}\,w^{n}{{\text{\ \ as }}}n\rightarrow\infty
	,\label{eq:asymp}%
	\end{equation}
	where $\theta_{2}$ and $w$ are given by (\ref{eq:theta2})--(\ref{eq:w}).
\end{theorem}

\begin{corollary}
	\label{c:asymp-monot}When $\theta_{1}\neq\theta_{0}$ the sequence
	${{{\mathbb{E}}}}_{\theta_{1}}[q_{n}^{\theta_{0}}]$ is eventually strictly decreasing.
\end{corollary}

\begin{proof}
	As $n\rightarrow\infty$ we have ${{{\mathbb{E}}}}_{\theta_{1}}[q_{n+1}%
	^{\theta_{0}}]/{{{\mathbb{E}}}}_{\theta_{1}}[q_{n}^{\theta_{0}}]\sim
	(\sqrt{n+1}/\sqrt{n})w\rightarrow w<1.$
\end{proof}

Note that the prior $\pi$ does not appear in (\ref{eq:asymp-theta0}); as in
the Bernstein--von Mises theorem (which is used in the proof below), the prior
is \textquotedblleft overwhelmed" by the data as the number of observations
increases; this is however \emph{not} the case in (\ref{eq:asymp}) when
$\theta_{1}\neq\theta_{0}$.

Before proving the theorem, we apply it to a number of classical examples:

\begin{itemize}
	\item \textbf{Bernoulli}$(\theta)$ for $\theta\in\Theta=(0,1).$ Here
	${{{\mathcal{X}}}}=\{0,1\},$ $\eta(\theta)=\log(\theta/(1-\theta)),$ and
	$T(x)=x,$ and so ${{\mathbb{V}}}ar_{\theta}(T(x))=\theta(1-\theta)$ and
	$\eta^{\prime}(\theta)=1/(\theta(1-\theta)),$ yielding $I(\theta
	)=1/(\theta(1-\theta))$ and
	\begin{align*}
	{{{\mathbb{E}}}}_{\theta_{0}}[q_{n}^{\theta_{0}}] &  \sim\frac{\sqrt{n}%
	}{2\sqrt{\uppi\,\theta_{0}(1-\theta_{0})}}{{{\text{\ \ \ (for }}}}0<\theta
	_{0}<1),\\
	{{{\mathbb{E}}}}_{\theta_{1}}[q_{n}^{\theta_{0}}] &  \sim\frac{\pi(\theta
		_{0})\sqrt{n}\,w^{n}}{2\pi(\theta_{2})\sqrt{\uppi\,\theta_{2}(1-\theta_{2})}%
	}{{{\text{\ \ (for }}}}\theta_{0}\neq\theta_{1}),
	\end{align*}
	where
	\[
	\theta_{2}=\frac{\sqrt{\theta_{0}\theta_{1}}}{\sqrt{\theta_{0}\theta_{1}%
		}+\sqrt{(1-\theta_{0})(1-\theta_{1})}}{{{\text{\ and\ \ }}}}w=\left(
	\sqrt{\theta_{0}\theta_{1}}+\sqrt{(1-\theta_{0})(1-\theta_{1})}\right)  ^{2}.
	\]

	\item \textbf{Normal}$(\theta,\sigma^{2})$ for $\theta\in\Theta={{{\mathbb{R}%
	}}}$ and fixed $\sigma>0.$ Here ${{{\mathcal{X}}}}={{{\mathbb{R}}}}$,
	$\eta(\theta)=\theta,$ and $T(x)=x/\sigma^{2},$ and so $I(\theta
	)={{\mathbb{V}}}ar_{\theta}(T(x))=1/\sigma^{2}$, yielding
	\begin{align*}
	{{{\mathbb{E}}}}_{\theta_{0}}[q_{n}^{\theta_{0}}] &  \sim\frac{\sqrt{n}%
	}{2\sigma\sqrt{\uppi}},\nonumber\\
	{{{\mathbb{E}}}}_{\theta_{1}}[q_{n}^{\theta_{0}}] &  \sim\frac{\pi(\theta
		_{0})\sqrt{n}}{2\pi(\theta_{2})\sigma\sqrt{\uppi}}\operatorname{exp}\left(
	-\frac{n(\theta_{0}-\theta_{1})^{2}}{4\sigma^{2}}\right)  ,
	\end{align*}
	where $\theta_{2}=(\theta_{0}+\theta_{1})/2.$
	
	\item \textbf{Exponential}$(\theta)$ for $\theta\in\Theta=(0,\infty)$ (i.e.,
	$p_{\theta}(x)=\theta e^{-\theta x}$ for $x>0)$. Here ${{{\mathcal{X}}}%
	}=(0,\infty),$ $\eta(\theta)=\theta,$ and $T(x)=-x,$ and so $I(\theta
	)={{\mathbb{V}}}ar_{\theta}(T(x))=1/\theta^{2}$, yielding
	\begin{align}\label{eq:expp}
	{{{\mathbb{E}}}}_{\theta_{0}}[q_{n}^{\theta_{0}}] &  \sim\frac{\sqrt{n}%
	}{2\theta_{0}\sqrt{\uppi}},\nonumber\\
	{{{\mathbb{E}}}}_{\theta_{1}}[q_{n}^{\theta_{0}}] &  \sim\frac{\pi(\theta
		_{0})\sqrt{n}}{2\pi(\theta_{2})\theta_{2}\sqrt{\uppi}}\left(  \frac{\theta
		_{0}\theta_{1}}{\theta_{2}^{2}}\right)  ^{n},
	\end{align}
	where $\theta_{2}=(\theta_{0}+\theta_{1})/2.$
	
	\item \textbf{Poisson}$(\theta)$ for $\theta\in\Theta=(0,\infty).$ Here
	${{{\mathcal{X}}}}={{{\mathbb{N}}}}$, $\eta(\theta)=\log\theta,$ and $T(x)=x,$
	and so ${{\mathbb{V}}}ar_{\theta}(T(x))=\theta$ and $\eta^{\prime}%
	(\theta)=1/\theta,$ yielding $I(\theta)=1/\theta$ and
	\begin{align*}
	{{{\mathbb{E}}}}_{\theta_{0}}[q_{n}^{\theta_{0}}]  &  \sim\frac{\sqrt{n}%
	}{2\sqrt{\uppi \,\theta_{0}}},\\
	{{{\mathbb{E}}}}_{\theta_{1}}[q_{n}^{\theta_{0}}]  &  \sim\frac{\pi(\theta
		_{0})\sqrt{n}}{2\pi(\theta_{2})\sqrt{\uppi \,\theta_{2}}}\operatorname{exp}%
	\left(  -n(\theta_{0}+\theta_{1}-2\theta_{2}\right)  ),
	\end{align*}
	where $\theta_{2}=\sqrt{\theta_{0}\theta_{1}}.$
\end{itemize}

\begin{proof}
	[Proof of Theorem \ref{th:asymp}]Part (ii) follows from part (i) by
	Proposition \ref{p:theta-bar} (when $\theta_{1}\neq\theta_{0}$ the point
	$\theta_{2}$ lies strictly between $\theta_{0}$ and $\theta_{1},$ and so is an
	interior point of $\Theta).$ We will thus prove (i). It is convenient to
	assume without loss of generality that $\eta(\theta)\equiv\theta$ (this is
	called the \textquotedblleft canonical" representation); indeed, since
	$\eta^{\prime}>0,$ the transformation $\tilde{\theta}:=\eta(\theta)$ (which
	preserves the convexity of the parameter space and maps interior points to
	interior points) yields: $\tilde{p}_{\tilde{\theta}}(x)=\operatorname{exp}%
	(\tilde{\theta}\cdot T(x)-A(\tilde{\theta})-B(x));\;$ $\tilde{\pi}%
	(\tilde{\theta})=\pi(\theta)/\eta^{\prime}(\theta);$\ $\tilde{q}_{n}%
	^{\tilde{\theta}}=q_{n}^{\theta}/\eta^{\prime}(\theta);$ and $\tilde{I}%
	(\tilde{\theta})=I(\theta)/(\eta^{\prime}(\theta))^{2};$ hence $\tilde{q}%
	_{n}^{\tilde{\theta}}/\sqrt{\tilde{I}(\tilde{\theta})}=q_{n}^{\theta}%
	/\sqrt{I(\theta)},$ and so (\ref{eq:asymp-theta0}) for $\theta$ is equivalent
	to (\ref{eq:asymp-theta0}) for $\tilde{\theta}.$ From now on we thus have
	\begin{equation}
	p_{\theta}(x)=\operatorname{exp}(\theta\,T(x)-A(\theta
	)-B(x))\label{eq:exp-fam-canon}%
	\end{equation}
	for all $x\in{{{\mathcal{X}}}}$ and $\theta\in\Theta,$ and so $I(\theta
	)=A^{\prime\prime}(\theta)$. For $s_{n}=(x_{1},...,x_{n})\in{{{\mathcal{X}}}%
	}^{n},$ let $\widehat{\theta}_{n}\equiv\widehat{\theta}_{n}(s_{n}):=\arg
	\max_{\theta\in\Theta}\sum_{i=1}^{n}\log p_{\theta}(x_{i})$ denote the maximum
	likelihood estimator (MLE); thus $\widehat{\theta}_{n}$ minimizes the strictly
	convex function $h_{n}(\theta):=A(\theta)-\theta\cdot\overline{t}_{n},$ where
	$\overline{t}_{n}\equiv\overline{t}_{n}(s_{n}):=(1/n)\sum_{i=1}^{n}T(x_{i});$
	if $\widehat{\theta}_{n}$ is an interior point of $\Theta$ then $h_{n}%
	^{\prime}(\widehat{\theta}_{n})=0,$ i.e., $A^{\prime}(\widehat{\theta}%
	_{n})=\overline{t}_{n}.$ Put ${{{\mathbb{P}}}}_{0}\equiv{{{\mathbb{P}}}%
	}_{\theta_{0}},$ $p_{0}\equiv p_{\theta_{0}},$ and ${{{\mathbb{E}}}}_{0}%
	\equiv{{{\mathbb{E}}}}_{\theta_{0}},$ respectively, for the probability, density, and
	expectation under $\theta_{0}$. Given iid observations $x_{1},x_{2},...,$
	under ${{\mathbb{P}}}_{0},$ we have (see, e.g., \cite{Schervish}, Theorem 7.57
	or \cite{FergLS}, Theorem 18)\footnote{Notation: $\overset{{{\cal{L}}%
		}}{\longrightarrow}$ means convergence in law (or distribution), and
		$\overset{{{{\mathbb{P}}}}_{0}}{\longrightarrow}$ means convergence in
		probability with respect to the probability ${{{\mathbb{P}}}}_{0}.$}
	$\sqrt{nI(\theta_{0})}(\widehat{\theta}_{n}-\theta_{0}%
	)\overset{{{\cal{L}}}}{\longrightarrow}$ $\mathcal{N}(0,1)$, which
	implies $\widehat{\theta}_{n}\overset{{{{\mathbb{P}}}}_{0}%
	}{\longrightarrow}\theta_{0}$ (the \textquotedblleft consistency" of the MLE). 
	
	For convenience we divide the proof of (\ref{eq:asymp-theta0}) (for
	$\theta_{0}\in{{{\mathrm{int}}}}\Theta)$ into a number of steps, as follows.
	First, we show that with high ${\mathbb{P}}_{0}$-probability the posterior
	$q_{n}^{\widehat{\theta}_{n}}$ at $\widehat{\theta}_{n}$ converges to the
	appropriate limit by the Bernstein--von Mises theorem (Step 1). Second, since
	$\widehat{\theta}_{n}-\theta_{0}$ is approximately normal (as seen above), we
	show that replacing $\widehat{\theta}_{n}$ with $\theta_{0}$ requires a factor
	of $1/\sqrt{2}$ on expectation (Steps 2 and 3). We then prove that the
	$\theta_{0}$-posterior $q_{n}^{\theta_{0}}$ is ${\mathrm{O}}(\sqrt{n})$ (Step
	4), and so sets of small ${\mathbb{P}}_{0}$-probability can be ignored (Step
	5); this completes the proof.
	
	\noindent$\bullet$ \emph{Step 1:} Let $J(\theta):=\sqrt{I(\theta)/(2\uppi)};$
	then
	\[
	\frac{1}{\sqrt{n}}q_{n}^{\widehat{\theta}_{n}}(s_{n})\overset{{{{\mathbb{P}}}%
		}_{0}}{\longrightarrow}J(\theta_{0}){{{\text{ as }}}}n\rightarrow\infty.
	\]

	\begin{proof}
		The Bernstein--von Mises theorem says that under $\theta_{0}$ the posterior
		density of $\vartheta:=\sqrt{nI(\widehat{\theta}_{n})}(\theta-\widehat{\theta
		}_{n})$ converges as $n\rightarrow\infty$ to the standard normal density
		$\varphi$; specifically, Theorem 7.89 in \cite{Schervish} (in Appendix
		\ref{sec-a:step1} we show that all the \textquotedblleft general regularity
		conditions" of the theorem hold) applied at $\vartheta=0,$ i.e.,
		$\theta=\widehat{\theta}_{n}$, yields $q_{n}^{\widehat{\theta}_{n}}%
		(s_{n})/\sqrt{nI(\widehat{\theta}_{n})}\overset{{{{\mathbb{P}}}}%
			_{0}}{\longrightarrow}\varphi(0)=1/\sqrt{2\uppi}$ . Now use $\widehat{\theta
		}_{n}\overset{{{{\mathbb{P}}}}_{0}}{\longrightarrow}\theta_{0}$ and the
		continuity of $I$.
	\end{proof}
	
	Given $0<\varepsilon<1,$ let $\Omega_{n}^{1}\equiv\Omega_{n}^{1}(\varepsilon)$
	be the event that $\left\vert (1/\sqrt{n})q_{n}^{\widehat{\theta}_{n}}%
	(s_{n})/J(\theta_{0})-1\right\vert \leq\varepsilon;$ thus ${{{\mathbb{P}}}%
	}_{0}(\Omega_{n}^{1})\rightarrow1$ (as $n\rightarrow\infty)$ by Step 1. Let
	$\Omega_{n}^{2}\equiv\Omega_{n}^{2}(\varepsilon)$ be the event that
	$\widehat{\theta}_{n}$ is an interior point of $\Theta$ and $|\pi(\theta
	_{0})/\pi(\widehat{\theta}_{n})-1|\leq\varepsilon;$ since $\widehat{\theta
	}_{n}\overset{{{{\mathbb{P}}}}_{0}}{\longrightarrow}\theta_{0}\in
	{{{\mathrm{int}}}}\Theta$ and $\pi$ is continuous and positive, ${{{\mathbb{P}%
	}}}_{0}(\Omega_{n}^{2})\rightarrow1.$ Put $\Omega_{n}:=\Omega_{n}^{1}%
	\cap\Omega_{n}^{2};$ then ${{{\mathbb{P}}}}_{0}(\Omega_{n})\rightarrow1.$
	
	\noindent$\bullet$ \emph{Step 2:}
	\[
	{{{\mathbb{E}}}}_{\theta_{0}}\left[  \frac{p_{0}(s_{n})}{p_{\widehat{\theta
			}_{n}}(s_{n})}{{{\mathbf{1}}}}_{\Omega_{n}}\right]  \rightarrow\frac{1}%
	{\sqrt{2}}{{{\text{ as }}}}n\rightarrow\infty.
	\]

	\begin{proof}
		Recall that $h_{n}(\theta):=A(\theta)-\theta\cdot\overline{t}_{n};$ then
		\[
		Y_{n}:=\frac{p_{0}(s_{n})}{p_{\widehat{\theta}_{n}}(s_{n})}{{{\mathbf{1}}}%
		}_{\Omega_{n}}=\operatorname{exp}\left(  -n[h(\theta_{0})-h(\widehat{\theta
		}_{n})]\right)  {{{\mathbf{1}}}}_{\Omega_{n}}.
		\]
		In $\Omega_{n}\subseteq\Omega_{n}^{2}$ the point $\widehat{\theta}_{n}$ is an
		interior point of $\Theta$ and so the Taylor series of $h$ around
		$\widehat{\theta}_{n},$ where $h^{\prime}(\widehat{\theta}_{n})=0$ and
		$h^{\prime\prime}(\widehat{\theta}_{n})=A^{\prime\prime}(\widehat{\theta}%
		_{n}),$ yields an intermediate point $\theta_{n}$ between $\theta_{0}$ and
		$\widehat{\theta}_{n}$ such that
		\[
		Y_{n}=\operatorname{exp}\left(  -\frac{1}{2}nA^{\prime\prime}(\theta
		_{n})(\theta_{0}-\widehat{\theta}_{n})^{2}\right)  {{{\mathbf{1}}}}%
		_{\Omega_{n}}.
		\]
		Now $\widehat{\theta}_{n}\overset{{{{\mathbb{P}}}}_{0}}{\longrightarrow}%
		\theta_{0}$ implies $\theta_{n}\overset{{{{\mathbb{P}}}}_{0}}{\longrightarrow
		}\theta_{0}$ and thus $A^{\prime\prime}(\theta_{n})\overset{{{{\mathbb{P}}}%
			}_{0}}{\longrightarrow}A^{\prime\prime}(\theta_{0})=I(\theta_{0});$ also,
		${{{\mathbf{1}}}}_{\Omega_{n}}\overset{{{{\mathbb{P}}}}_{0}}{\longrightarrow
		}1$ (because ${{{\mathbb{P}}}}_{0}(\Omega_{n})\rightarrow1).$ Under $p_{0}$ we
		have $\sqrt{nI(\theta_{0})}(\widehat{\theta}_{n}-\theta_{0}%
		)\overset{{{\cal{L}}}}{\longrightarrow}$ $Z\equiv\mathcal{N}(0,1),$ as
		mentioned before Step 1; altogether, $Y_{n}\overset{{{\cal{L}}%
		}}{\longrightarrow}\operatorname{exp}(-Z^{2}/2)$. The $Y_{n}$ are uniformly
		bounded ($0\leq Y_{n}\leq1,$ because $A^{\prime\prime}>0),$ and so
		${{{\mathbb{E}}}}_{0}[Y_{n}]\rightarrow{{{\mathbb{E[}}}}\operatorname{exp}%
		(-Z^{2}/2)]=1/\sqrt{2}.$
	\end{proof}
	
	\noindent$\bullet$ \emph{Step 3:}
	\[
	\limsup_{n\rightarrow\infty}\left\vert \frac{1}{\sqrt{n}}{{{\mathbb{E}}}%
	}_{\theta_{0}}[q_{n}^{\theta_{0}}{{{\mathbf{1}}}}_{\Omega_{n}}]-\frac
	{J(\theta_{0})}{\sqrt{2}}\right\vert \leq\varepsilon^{\prime},
	\]
	where $\varepsilon^{\prime}:=(3J(\theta_{0})/\sqrt{2})\varepsilon.$
	
	\begin{proof}
		For every $s_{n}$ we have
		\[
		\frac{1}{\sqrt{n}}q_{n}^{\theta_{0}}(s_{n})=\frac{q_{n}^{\widehat{\theta}_{n}%
			}(s_{n})}{\sqrt{n}}\cdot\frac{\pi(\theta_{0})}{\pi(\widehat{\theta}_{n})}%
		\cdot\frac{p_{0}(s_{n})}{p_{\widehat{\theta}_{n}}(s_{n})}.
		\]
		In $\Omega_{n}^{1}$ the first factor is at most $(1+\varepsilon)J(\theta
		_{0}),$ and in $\Omega_{n}^{2}$ the second factor is at most $1+\varepsilon$;
		hence
		\[
		\frac{1}{\sqrt{n}}q_{n}^{\theta_{0}}(s_{n}){{{\mathbf{1}}}}_{\Omega_{n}}%
		\leq(1+\varepsilon)^{2}J(\theta_{0})\frac{p_{0}(s_{n})}{p_{\widehat{\theta
				}_{n}}(s_{n})}{{{\mathbf{1}}}}_{\Omega_{n}}.
		\]
		Taking expectation under $\mathbb{P}_{0}$ yields by Step 2
		\[
		\limsup_{n\rightarrow\infty}\frac{1}{\sqrt{n}}{{{\mathbb{E}}}}_{\theta_{0}%
		}[q_{n}^{\theta_{0}}{{{\mathbf{1}}}}_{\Omega_{n}}]\leq(1+\varepsilon)^{2}%
		\frac{J(\theta_{0})}{\sqrt{2}}\leq\frac{J(\theta_{0})}{\sqrt{2}}%
		+\varepsilon^{\prime}.
		\]
		Similarly,
		\[
		\liminf_{n\rightarrow\infty}\frac{1}{\sqrt{n}}{{{\mathbb{E}}}}_{\theta_{0}%
		}[q_{n}^{\theta_{0}}{{{\mathbf{1}}}}_{\Omega_{n}}]\geq(1-\varepsilon)^{2}%
		\frac{J(\theta_{0})}{\sqrt{2}}\geq\frac{J(\theta_{0})}{\sqrt{2}}%
		-\varepsilon^{\prime},
		\]
		completing the proof.
	\end{proof}
	
	\noindent$\bullet$ \emph{Step 4:} There is a constant $C<\infty$ such that
	\[
	q_{n}^{\theta_{0}}(s_{n})\leq C\sqrt{n}
	\]
	for all $s_{n}$ and all $n\geq1.$
	
	\begin{proof}
		We will show that there is a constant $c>0$ such that
		\[
		\frac{\pi(\theta_{0})}{q_{n}^{\theta_{0}}(s_{n})}=\frac{p(s_{n})}{p_{0}%
			(s_{n})}\geq\frac{c}{\sqrt{n}}%
		\]
		for all $n\geq1.$ Take $\delta>0$ such that $[\theta_{0}-\delta,\theta
		_{0}+\delta]\subset{{{\mathrm{int}}}}\Theta;$ then
		\[
		\frac{p(s_{n})}{p_{0}(s_{n})}=\int_{\Theta}\operatorname{exp}(-n[h_{n}%
		(\theta)-h_{n}(\theta_{0})])\pi(\theta)\,{{{\mathrm{d}}}}\theta\geq\rho
		\int_{\theta_{0}-\delta}^{\theta_{0}+\delta}H_{n}(\theta)\,{{{\mathrm{d}}}%
		}\theta,
		\]
		where $H_{n}(\theta):=\operatorname{exp}(-n[h_{n}(\theta)-h_{n}(\theta_{0})])$
		and $\rho:=\min_{\theta\in\lbrack\theta_{0}-\delta,\theta_{0}+\delta]}%
		\pi(\theta)>0$ (recall that the density $\pi$ is continuous and strictly
		positive on $\Theta).$ The strictly convex function $h_{n}(\theta
		)=A(\theta)-\theta\cdot\overline{t}_{n}$ has a unique minimizer in
		$[\theta_{0}-\delta,\theta_{0}+\delta],$ call it $\xi;$ without loss of
		generality assume that $\xi\geq\theta_{0}.$ When $\xi\geq\theta_{0}+\delta$
		the function $h_{n}$ is decreasing for $\theta\leq\xi,$ and thus for all
		$\theta\in\lbrack\theta_{0},\theta_{0}+\delta]$ we have $h_{n}(\theta)\leq
		h_{n}(\theta_{0}),$ and hence $H_{n}(\theta)\geq1,$ which gives $p(s_{n}%
		)/p_{0}(s_{n})\geq\rho\delta\geq\rho\delta/\sqrt{n}$ for all $n\geq1.$ When
		$\theta_{0}\leq\xi<\theta_{0}+\delta$ we have\footnote{In this case, where
			$\xi$ is an interior point, $\xi$ is the minimum over all $\Theta,$ and thus
			$\xi=\widehat{\theta}_{n}.$ The argument here is an instance of the Laplace
			method.} $h_{n}^{\prime}(\xi)=0$ and so $h_{n}(\theta)-h_{n}(\theta_{0})\leq
		h_{n}(\theta)-h_{n}(\xi)=h_{n}^{\prime\prime}(\zeta)(\theta-\xi)^{2}/2$ for
		every $\theta\in\lbrack\theta_{0}-\delta,\theta_{0}+\delta]$ (by the
		second-order Taylor expansion), where $\zeta\equiv\zeta_{\theta}$ is an
		intermediate point between $\theta$ and $\xi,$ and so $\zeta\in\lbrack
		\theta_{0}-\delta,\theta_{0}+\delta].$ Since $h_{n}^{\prime\prime}%
		=A^{\prime\prime}$ we get $0<h_{n}^{\prime\prime}(\zeta)\leq\alpha
		:=\max_{\theta\in\lbrack\theta_{0}-\delta,\theta_{0}+\delta]}A^{\prime\prime
		}(\theta),$ and so $h_{n}(\theta)-h_{n}(\theta_{0})\leq\alpha(\theta-\xi
		)^{2}/2$ and
		\[
		\rho\int_{\theta_{0}-\delta}^{\theta_{0}+\delta}H(\theta)\,{{{\mathrm{d}}}%
		}\theta\geq\rho\int_{\xi-\delta}^{\xi}\operatorname{exp}\left(  -\frac
		{n\alpha}{2}(\theta-\xi)^{2}\right)  {{\mathrm{d}}}\theta=\frac{\rho
			\sqrt{2\uppi}}{\sqrt{n\alpha}}\left(  \frac{1}{2}-\Phi(-\delta\sqrt{n\alpha
		})\right)
		\]
		(because $[\xi-\delta,\xi]\subset\lbrack\theta_{0}-\delta,\theta_{0}+\delta];$
		here $\Phi$ denotes the cumulative standard normal distribution). Since
		$\Phi(-\delta\sqrt{n\alpha})\rightarrow0$ as $n\rightarrow\infty,$ the final
		expression is $\sim c/\sqrt{n}$ for some $c>0,$ which completes the proof.
	\end{proof}
	
	\noindent$\bullet$ \emph{Step 5:} Let $\Omega_{n}^{c}$ denote the complement of $\Omega_{n}$. We have
	\[
	\frac{1}{\sqrt{n}}{{{\mathbb{E}}}}_{\theta_{0}}[q_{n}^{\theta_{0}%
	}{{{\mathbf{1}}}}_{\Omega_{n}^{c}}]\rightarrow0{{{\text{ \ as }}}}%
	n\rightarrow\infty.
	\]

	\begin{proof}
		Use the uniform boundedness of $(1/\sqrt{n})q_{n}^{\theta_{0}}$ by Step 4 and
		${{{\mathbb{P}}}}_{0}(\Omega_{n})\rightarrow1$.
	\end{proof}
	
	Adding the results of Steps 3 and 5 and noting that $\varepsilon>0$ is
	arbitrary yields (\ref{eq:asymp-theta0}), and thus completes the proof of Theorem \ref{th:asymp}.
\end{proof}

See Appendix \ref{sec-a:asymp} for additional comments, extensions, and
technical details.

\section{Log-concavity}

\label{sec:loggcon} In this section we discuss setups in which the sequence
${{\mathbb{E}}}_{\theta_{1}}[q_{n}^{\theta_{0}}]$ is log-concave in $n$. A
sequence of positive numbers $\psi(n)$ for $n\geq1$ is (\emph{strictly}%
)\emph{\ log-concave} if $\log\psi(n)$ is a (strictly) concave function of
$n$; this is equivalent to $\psi(n)^{2}\geq\psi(n-1)\psi(n+1)$ for every
$n\geq2$, with $>$ for the strict version. The sequence $\psi(n)$
is\emph{\ unimodal} if there exists $0\leq n_{0}\leq\infty$ (possibly equal to
$0$ or $\infty)$ such that $\psi(n)$ is increasing for $n\leq n_{0}$ and
decreasing for $n\geq n_{0}.$ Log-concavity clearly implies unimodality (with
$n_{0}$ that maximizes $\psi(n)$).

In Corollary \ref{c:asymp-monot} we saw that for $\theta_{1}\neq\theta_{0}$
the sequence $\psi_{\theta_{0},\theta_{1}}(n)={{\mathbb{E}}}_{\theta_{1}%
}[q_{n}^{\theta_{0}}]$ is eventually strictly decreasing. We will now show
that in certain natural setups, with conjugate priors, this can be
strengthened to unimodality, and, in fact, to log-concavity, with respect to
the time period $n.$ We do this in three setups. The first one consists of iid
Bernoulli observations with a uniform prior; the second, of iid normal
observations with a normal prior; and the third, of iid exponential observations
with an exponential prior. In the normal case log-concavity is obtained only
from some $n_{0}$ on, where $n_{0}$ may be $1$ or arbitrarily large, depending
on the parameters; Figure 3 in Section \ref{sec:nonmonmult} is typical of the
latter case. The analysis suggests that general log-concavity results may be
hard to obtain, as indicated by the different proofs in the three cases, as
well as by the dependence of the result on the specific prior (see the normal
case, or take a Beta prior in the Bernoulli case).

In Proposition \ref{p:theta-bar} of Section \ref{sec:asympt} we get a
log-linear relation between $\psi_{\theta_{0},\theta_{1}}(n)$ and
$\psi_{\theta_{2},\theta_{2}}(n),$ and so it suffices to consider only the
case where $\theta_{1}=\theta_{0}.$

\begin{corollary}
	\label{cor:suffsame} Under the assumptions of Proposition \ref{p:theta-bar},
	the sequence $\psi_{\theta_{0},\theta_{1}}(n)$ is log-concave/convex in $n$ if
	and only if the sequence $\psi_{\theta_{2},\theta_{2}}(n)$ is
	log-concave/convex in $n.$
\end{corollary}

\subsection{Bernoulli observations with uniform prior}

\label{sec:uniprior}

This section deals with sequences of iid Bernoulli observations with a
parameter $\theta$ that is uniformly distributed in $(0,1).$ We show that in
this case $\psi_{\theta_{0},\theta_{1}}(n)={{\mathbb{E}}}_{\theta_{1}}%
[q_{n}^{\theta_{0}}]$ is a log-concave function of the time period $n$, and so
unimodal in $n.$ We obtain this result by proving in Section \ref{sec:revtur}
a \textquotedblleft reversal" of the reverse Tur\'{a}n inequality for Legendre
polynomials, which may be of independent interest.

\begin{theorem}
	\label{th:binunprlog} Let $x_{1},x_{2},\ldots$ be iid ${{\mathrm{Bernoulli}}%
	}(\theta)$ observations, and let the prior distribution of $\theta$ be the
	uniform distribution on $\Theta=(0,1).$ For every $\theta_{0}$ and $\theta
	_{1}$ in $\Theta$ the sequence $\psi_{\theta_{0},\theta_{1}}(n)={{\mathbb{E}}%
	}_{\theta_{1}}[q_{n}^{\theta_{0}}]$ is log-concave for $n\geq1$ \textbf{(}and
	strictly log-concave for $n\geq2$\textbf{)}, and hence unimodal.
\end{theorem}

\begin{proof}
	As in Section \ref{sec:nonmonmult}, we work with the sufficient statistic
	$u_{n}:=\sum_{i=1}^{n}x_{i}$, whose distribution given $\theta$ is
	\textrm{Binomial}$(n,\theta).$ The marginal distribution of $u_{n}$ when the
	prior is uniform on $(0,1)$ is then the uniform distribution on the set
	$\{0,1,...,n\};$ i.e., ${{\mathbb{P}}}(u_{n}=k)=1/(n+1)$ for $0\leq k\leq n$
	(this is a well-known result; see, e.g., \cite{GG}, page 155). Therefore
	\[
	\psi_{\theta,\theta}(n)=\sum_{k=0}^{n}\frac{\left(  {{\mathbb{P}}}_{\theta
		}(u_{n}=k)\right)  ^{2}}{{{\mathbb{P}}}(u_{n}=k)}=(n+1)\sum_{k=0}^{n}%
	{\binom{n}{k}}^{2}\theta^{2k}(1-\theta)^{2(n-k)}.
	\]
	The log-concavity of the last expression follows from Corollary \ref{c:PSI} in
	the next section, with $y=\theta^{2}$ and $z=(1-\theta)^{2}.$ Finally, for
	$\theta_{1}\neq\theta_{0}$ apply Corollary \ref{cor:suffsame}.
\end{proof}

It is natural to try to generalize from the uniform prior to other priors,
such as Beta distributions (the class of conjugate priors here). However,
numerical calculations show that ${{\mathbb{E}}}_{\theta_{1}}[q_{n}%
^{\theta_{0}}]$ need \emph{not} be log-concave: take, for example, the
${{\mathrm{Beta}}}(7,1)$ prior, $\theta_{0}=3/4,$ $\theta_{1}=9/10,$ and
$n=2,3,4$ (this is clearly robust to small changes in the parameters).

\subsection{Reversing the reverse Tur\'{a}n inequality for Legendre
	polynomials}

\label{sec:revtur} This section proves an interesting reversal of the reverse
Tur\'{a}n inequality for Legendre polynomials for $|x|>1;$ it yields in
particular the log-concavity of the previous section.

The Legendre polynomial of degree $n\geq0$ is defined as
\[
P_{n}(x){\;:=\;}{\frac{1}{2^{n}}}\sum_{k=0}^{n}{\binom{n}{k}}^{2}%
(x-1)^{n-k}(x+1)^{k};
\]
see, e.g., \cite{szw} and the references therein. The well-known Tur\'{a}n
inequality for Legendre polynomials, first published in \cite{sze}, states
that
\begin{equation}
P_{n}^{2}(x)\geq P_{n-1}(x)P_{n+1}(x)\,\,{{\text{ for all }}}\,|x|\leq1
\label{eq:turan}%
\end{equation}
holds for every $n\geq1,$ with equality if and only if\footnote{Put $0^{0}=1$;
	then $P_{n}(1)=1$ and $P_{n}(-1)=(-1)^{n}.$} $|x|=1$. Its reverse version,
\begin{equation}
P_{n}^{2}(x)<P_{n-1}(x)P_{n+1}(x)\,\,{{\text{ for all }}}\,|x|>1,
\label{eq:anti-turan}%
\end{equation}
holds for every $n\ge1;$ see, e.g., \cite{szw}, Theorem 1.

We next show that multiplying $P_{n}$ by $n+1$ reverses (\ref{eq:anti-turan})
for all\footnote{Starting with $P_{1}$ rather than $P_{0}$ (which is what is
	needed for our Theorem \ref{th:binunprlog}).} $|x|>1$.

\begin{theorem}
	\label{th:turan} The inequality
	\begin{equation}
	\frac{P_{n-1}(x)P_{n+1}(x)}{P_{n}^{2}(x)}\leq\frac{(n+1)^{2}}{n(n+2)}%
	\,\,{{\text{ for all }}}\,|x|>1 \label{eq:beta}%
	\end{equation}
	holds for every $n\geq2,$ with equality if and only if $n=2$ and $|x|=\sqrt
	{3}$.
\end{theorem}

\begin{proof}
	Putting
	\[
	R_{n}{\;:=\;}\frac{P_{n-1}P_{n+1}}{P_{n}^{2}}\;\;{{\text{and\ \ }}}%
	a_{n}{\;:=\;}\frac{(n+1)^{2}}{n(n+2)}
	\]
	we will prove that
	\[
	1<R_{n}(x)\leq a_{n}
	\]
	for all $|x|>1$ and $n\geq2;$ the first inequality (which also holds for
	$n=1)$ is the reverse Tur\'{a}n inequality (\ref{eq:anti-turan}), for which we
	provide a simple proof as well. Because $P_{n}(-x)=(-1)^{n}P_{n}(x)$ it
	suffices to consider the range $x>1,$ where $P_{n}(x)>0$ for all $n\geq1$. Let
	$b_{n}$ be the leading coefficient of $P_{n}$, i.e., the coefficient of its
	highest power, $x^{n}$; then
	\[
	b_{n}=\frac{1}{2^{n}}\sum_{k=0}^{n}\binom{n}{k}^{2}=\frac{1}{2^{n}}\binom
	{2n}{n},
	\]
	and thus
	\begin{equation}
	R_{n}(\infty){\;:=\;}\lim_{x\rightarrow\infty}R_{n}(x)=\frac{b_{n-1}b_{n+1}%
	}{b_{n}^{2}}=\frac{n(2n+1)}{(n+1)(2n-1)}. \label{eq:nnn}%
	\end{equation}
	It is straightforward to verify that
	\begin{equation}
	1<R_{n}(\infty)<a_{n} \label{eq:limit}%
	\end{equation}
	for every $n\geq2$ (the first inequality is equivalent to $(2n^{2}%
	+n)/(2n^{2}+n-1)>1,$ and the second to $n^{2}-n-1>0).$ Next, differentiating
	$\log R_{n}(x)$ with respect to $x$ yields
	\[
	(\log R_{n})^{\prime}=\frac{R_{n}^{\prime}}{R_{n}}=\frac{P_{n-1}^{\prime}%
	}{P_{n-1}}+\frac{P_{n+1}^{\prime}}{P_{n+1}}-\frac{2P_{n}^{\prime}}{P_{n}}.
	\]
	Using  the following well-known formula (see, e.g., \cite{Ar}, Equation (12.26)
	for a convenient reference)
	\[
	(x^{2}-1)P_{n}^{\prime}(x)=nxP_{n}(x)-nP_{n-1}(x)
	\]
	for every $n\geq1$ and every $x$ we then obtain
	\begin{multline}
	(x^{2}-1)\frac{R_{n}^{\prime}}{R_{n}}%\\
	=(n-1)x-(n-1)\frac{P_{n-2}}{P_{n-1}}+(n+1)x-(n+1)\frac{P_{n}}{P_{n+1}%
	}-2nx+2n\frac{P_{n-1}}{P_{n}}\\
	=2n\frac{P_{n-1}}{P_{n}}-(n-1)\frac{P_{n-2}}{P_{n-1}}-(n+1)\frac{P_{n}%
	}{P_{n+1}}
	=(n+1)\frac{P_{n-1}}{P_{n}}\left(  \frac{2n}{n+1}-\frac{n-1}{n+1}R_{n-1}%
	-\frac{1}{R_{n}}\right)  . \label{eq:diff=0}%
	\end{multline}
	The proof of \eqref {eq:beta} is by induction on $n,$ separately for each one
	of the two inequalities. For the first inequality, assume by induction that
	$R_{n-1}(x)>1$ for every $x>1.$ If $R_{n}(x)\leq1$ for some $x>1,$ then, since
	$R_{n}(1)=1$ and $R_{n}(\infty)>1$ (see \eqref {eq:nnn}), there is $x_{\ast
	}>1$ where $R_{n}$ attains its minimum, and so $R_{n}(x_{\ast})\leq1$ and
	$R_{n}^{\prime}(x_{\ast})=0.$ Using (\ref{eq:diff=0}) and then $R_{n-1}%
	(x_{\ast})>1$ (by the induction hypothesis) yields
	\[
	\frac{1}{R_{n}(x_{\ast})}=\frac{2n}{n+1}-\frac{n-1}{n+1}R_{n-1}(x_{\ast
	})<\frac{2n}{n+1}-\frac{n-1}{n+1}=1,
	\]
	in contradiction to $R_{n}(x_{\ast})\leq1.$ The induction starts with $n=1,$
	where we have
	\[
	P_{0}P_{2}-P_{1}^{2}=1\cdot\frac{1}{2}(3x^{2}-1)-\left(  x\right)  ^{2}%
	=\frac{1}{2}(x^{2}-1)>0{{\text{ for }}}x>1,
	\]
	and so $R_{1}(x)>1$ for every $x>1.$ For the second inequality, assume by
	induction that $R_{n-1}(x)<a_{n-1}$ for every $x>1.$ If $R_{n}(x)\geq a_{n}$
	for some $x>1,$ then, because $R_{n}(\infty)<a_{n}$ (see \eqref {eq:nnn} and
	(\ref{eq:limit})), there is $x^{\ast}>1$ where $R_{n}$ attains its maximum,
	and so $R_{n}(x^{\ast})\geq a_{n}$ and $R_{n}^{\prime}(x^{\ast})=0.$ Using
	\eqref {eq:diff=0} whose left-hand side vanishes at $x^{\ast}$, and then
	$R_{n-1}(x^{\ast})<a_{n-1}$ (by the induction hypothesis) yields
	\begin{align*}
	\frac{1}{R_{n}(x^{\ast})}  &  =\frac{2n}{n+1}-\frac{n-1}{n+1}R_{n-1}(x^{\ast
	})>\frac{2n}{n+1}-\frac{n-1}{n+1}a_{n-1}\\
	&  =\frac{2n}{n+1}-\frac{n-1}{n+1}\frac{n^{2}}{(n-1)(n+1)}=\frac
	{n(n+2)}{(n+1)^{2}}=\frac{1}{a_{n}},
	\end{align*}
	in contradiction to $R_{n}(x^{\ast})\geq a_{n}.$ The induction now starts with
	$n=2,$ where we have
	\[
	P_{1}P_{3}-\frac{9}{8}P_{2}^{2}=x\cdot\frac{1}{2}(5x^{3}-3x)-\frac{9}%
	{8}\left(  \frac{1}{2}(3x^{2}-1)\right)  ^{2}=-\frac{(x^{2}-3)^{2}}{32}%
	\leq0{{\text{ for }}}x>1,
	\]
	and so $R_{2}(x)\leq a_{2}$ for every $x>1,$ with equality only for
	$x=\sqrt{3};$ for $n=3$ we have $R_{3}(\sqrt{3})=19/18<16/15=a_{3},$ and so
	$x^{\ast}$ cannot be $\sqrt{3},$ and thus $R_{2}(x^{\ast})<a_{2}$ and the
	induction argument above gives $R_{3}(x)<a_{3}$ for every $x>1,$ and then
	$R_{n}(x)<a_{n}$ for every $x>1$ and $n\geq4.$
\end{proof}

\begin{corollary}
	\label{c:PHI}For fixed $x>1,$ the sequence $Q_{n}(x){\;:=\;}(n+1)P_{n}(x)$ is
	log-concave in $n$ for $n\geq1,$ and strictly log-concave for $n\geq2.$
\end{corollary}

Thus, $Q_{n}^{2}(x)\geq Q_{n-1}(x)Q_{n+1}(x)$ holds for every $n\geq2,$ with
strict inequality for $n\geq3.$ The result holds for $x=1$ as well, where
$P_{n}(1)=1$ for all $n.$

The sequence $P_{n}(x),$ for fixed $x\geq1,$ is log-convex in $n$ by the
reverse Tur\'{a}n inequality (\ref{eq:anti-turan}); Theorem \ref{th:turan}
says that multiplying $P_{n}(x)$ by $n+1$ makes the sequence log-concave
instead. Moreover, $(n+1)P_{n}(x)$ is the \textquotedblleft right" multiple
for which this reversal occurs: taking any smaller multiple, such as
$nP_{n}(x),$ also reverses the inequality (by (\ref{eq:beta}) and
$\frac{(n+1)^{2}}{n(n+2)}\leq\frac{n^{2}}{(n-1)(n+1)}$), whereas any larger
multiple, such as $(n+2)P_{n}(x),$ does not (consider $n=2$ and $x=\sqrt{3}).$

\begin{corollary}
	\label{c:PSI}For fixed $y,z\geq0$ and not both $0,$ the sequence
	\[
	S_{n}(y,z){\;:=\;}(n+1)\sum_{k=0}^{n}{\binom{n}{k}}^{2}y^{k}z^{n-k}
	\]
	is log-concave in $n$ for $n\geq1,$ and strictly log-concave for $n\geq2$ .
\end{corollary}

\begin{proof}
	When $y\neq z,$ say $y>z,$ put $x=(y+z)/(y-z),$ and then $x\geq1$ and
	$S_{n}(y,z)=(y-z)^{n}(n+1)P_{n}\left(  x\right)  =(y-z)^{n}Q_{n}(x)$, and we
	use Corollary \ref{c:PHI} for $x>1,$ and $P_{n}(1)=1$ for $x=1.$ When $y=z>0$
	we have $S_{n}(y,y)=y^{n}(n+1)\binom{2n}{n},$ and then $S_{n}^{2}%
	>S_{n-1}S_{n+1}$ is obtained from (\ref{eq:nnn})--(\ref{eq:limit}) or by
	direct calculation.
\end{proof}

\subsection{Normal observations with normal prior}

\label{sec:normnorm}

We now consider normal observations whose mean is normally distributed;
specifically, $x_{1},x_{2},\ldots$ are iid ${{\mathcal{N}}}(\theta,\sigma
^{2})$ observations (with $\sigma>0$ fixed), and the prior on $\theta
{{\mathbb{\ }}}$is the standard normal distribution ${{\mathcal{N}}}(0,1)$. We
show that the expected posterior $\psi_{\theta_{0},\theta_{1}}(n)\equiv
{{\mathbb{E}}}_{\theta_{1}}[q_{n}^{\theta_{0}}]$ is either a log-concave
function of $n,$ or a log-convex function up to some point and a log-concave
function thereafter (as in Figure 3); which case it is depends on $\theta
_{0},\theta_{1},$ and $\sigma$.

\begin{theorem}
	\label{th:normal} Let $x_{1},x_{2},\ldots$ be iid ${{\mathcal{N}}}%
	(\theta,\sigma^{2})$ observations (with $\sigma>0$ fixed), and let the prior
	distribution of $\theta$ be the normal standard distribution ${{\mathcal{N}}%
	}(0,1)$ on $\Theta={{\mathbb{R}}}.$
	
	(i) For every real $\theta$ and $n\geq1$ we have
	\begin{equation}
	\psi_{\theta,\theta}(n)=\frac{(n+\sigma^{2})}{\sigma\sqrt{2\uppi \,(2n+\sigma
			^{2})}}\operatorname{exp}\left(  -\frac{\theta^{2}\sigma^{2}}{4n+2{\sigma}%
		^{2}}\right)  , \label{eq:normal}%
	\end{equation}
	and thus
	\begin{equation}
	\psi_{\theta_{0},\theta_{1}}(n)=\operatorname{exp}\left(  \frac{\theta
		^{2}-\theta_{0}^{2}}{2}\right)  \psi_{\theta,\theta}(n)\,{\operatorname{exp}%
	}\left(  -\frac{n(\theta_{0}-\theta_{1})^{2}}{4\sigma^{2}}\right)
	\label{eq:nicenorm}%
	\end{equation}
	for every real $\theta_{0}$ and $\theta_{1}$ with $(\theta_{0}+\theta
	_{1})/2=\theta.$
	
	(ii) There is $n_{0}\geq0$ that depends on $\sigma$ and $\theta$ such that the
	sequence $\psi_{\theta,\theta}(n)$ is strictly log-convex for $n<n_{0}$ and
	strictly log-concave for $n>n_{0},$ and thus so are the sequences
	$\psi_{\theta_{0},\theta_{1}}(n)$ for every $\theta_{0},\theta_{1}$ with
	$(\theta_{0}+\theta_{1})/2=\theta.$
	
	(iii) When $\sigma^{2}\leq\sqrt{2}$ or $|\theta|\geq1/2$ the sequence
	$\psi_{\theta,\theta}(n)$ is strictly log-concave for $n\geq1.$
	
	(iv) For every $\theta_{1}\neq\theta_{0}$ the sequence $\psi_{\theta
		_{0},\theta_{1}}(n)$ has at most two critical points, and so it is of one of
	three types: always decreasing; increasing and then decreasing; decreasing,
	increasing, and then decreasing.
\end{theorem}

\begin{proof}
	(i) Take the sufficient statistic $u_{n}:=\sum_{i=1}^{n}x_{i}.$ The
	distribution of $u_{n}$ given $\theta$ is ${{\mathcal{N}}}(n\theta,n\sigma
	^{2}),$ and the marginal distribution of $u_{n}$ is ${{\mathcal{N}}}%
	(0,n^{2}+n\sigma^{2})$ (express $u_{n}$ as the sum of $n\theta$ and
	$u_{n}-n\theta).$ The result (\ref{eq:normal}) is then obtained by a standard
	computation, which we relegate to Appendix \ref{sec-a:normnorm}; as for
	(\ref{eq:nicenorm}), it then follows from Proposition \ref{p:theta-bar}.
	
	(ii) Let $\xi(n):=\log\psi_{\theta,\theta}(n).$ Taking derivatives with
	respect to $n$ (which we view as a continuous variable in (\ref{eq:normal}))
	yields $\xi^{\prime\prime}(n)=\gamma(n)/(2n+\sigma^{2})^{3},$ where
	\[
	\gamma(n)={\frac{({\sigma}^{4}-2n^{2})\left(  2n+{\sigma}^{2}\right)
		}{\left(  n+{\sigma}^{2}\right)  ^{2}}}-4\,{\theta}^{2}{\sigma}^{2}.
	\]
	The function $\gamma(n)$ is strictly decreasing for $n>0$ (because
	$\gamma^{\prime}(n)=-2n(2n^{2}+6n\sigma^{2}+3\sigma^{4})/(n+\sigma^{2}%
	)^{3}<0),$ and is negative for large enough $n$ (for sure when $2n^{2}%
	\geq\sigma^{4}).$ Thus either $\xi^{\prime\prime}$ is always negative, or it
	changes sign once from positive to negative, which means that either $\xi$ is
	always concave, or it is first convex and then concave. 
	
	(iii) If $\sigma^{2}\leq\sqrt{2}$ then $\sigma^{4}\leq2n^{2}$ for all
	$n\geq1,$ and so $\gamma(n)<0$ for all $n\geq1.$ If $|\theta|\geq1/2$ then
	$\gamma(0)\leq\sigma^{2}-4\theta^{2}\sigma^{2}\leq0,$ and so $\gamma(n)<0$ for
	all $n>0.$ 
	
	(iv) Let $\tilde{\xi}(n):=\log\psi_{\theta_{0},\theta_{1}}(n)$ and $\xi(n):=$
	$\log\psi_{\theta,\theta}(n)$ for $\theta=(\theta_{0}+\theta_{1})/2;$ then
	$\tilde{\xi}^{\prime}=\xi^{\prime}+\log w$ (for the appropriate $w)$ and
	$\tilde{\xi}^{\prime\prime}=\xi^{\prime\prime}.$ Since $\tilde{\xi}%
	^{\prime\prime}=\xi^{\prime\prime}$ changes sign at most once by (ii),
	$\tilde{\xi}^{\prime}$ can vanish at most twice, and so $\tilde{\xi}$ has at
	most two critical points; the last one must be a maximum since $\tilde{\xi
	}(n)\rightarrow-\infty$ as $n\rightarrow\infty.$ The same then holds for
	$\psi_{\theta_{0},\theta_{1}}=\exp\tilde{\xi}.$
\end{proof}

Figure 3 provides an example with two critical points, a minimum followed by a
maximum, which is thus the most that one may get in this normal setup.

\textbf{General normal prior: }When the prior on $\theta$ is a general
${{\mathcal{N}}}(\mu,\sigma_{\Theta}^{2})$ distribution, and the observations
$x_{i}$ given $\theta$ are iid ${{\mathcal{N}}}(\theta,\sigma_{X}^{2})$ (where
$\mu,\sigma_{\Theta}^{2},\sigma_{X}^{2}$ are fixed), the linear transformation
$x\rightarrow(x-\mu)/\sigma_{\Theta}$ reduces it to the above case with $\sigma
^{2}=\sigma_{X}^{2}/\sigma_{\Theta}^{2}$. Thus, $\sigma^{2}$ is now the ratio
of the variance of the observations to the variance of the prior. The above result
(iii), for instance, says that when the variance of the observations is
not too large relative to the variance of the prior (specifically, when
$\sigma_{X}^{2}\leq\sqrt{2}\sigma_{\Theta}^{2}),$ the sequence of expected
posteriors is log-concave, and thus unimodal, for $n\geq1$.

\subsection{Exponential observations with exponential prior\label{sec:expexp}}

In this section we consider exponential observations whose parameter is also
exponentially distributed; specifically, $x_{1},x_{2},\ldots$ are iid
${{\mathrm{Exp}}}(\theta)$ observations, and the prior distribution of
$\theta{{\mathbb{\ }}}$is ${{\mathrm{Exp}}}(1).$ Thus $p_{\theta}%
(x)=\theta\operatorname{exp}(-\theta x)$ and $\pi(\theta)=\operatorname{exp}%
(-\theta)$ for $x\geq0$ and $\theta>0$ (with $\Theta=(0,\infty)).$ Working
again with the sufficient statistic $u_{n}:=\sum_{i=1}^{n}x_{i},$ whose
distribution conditional on $\theta$ is the $\Gamma(n,\theta)$ distribution
(the $n$-fold convolution of ${{\mathrm{Exp}}}(\theta)\equiv\Gamma
(1,\theta)),$ we have $p_{\theta}(u_{n})=\theta^{n}u_{n}^{n-1}%
\operatorname{exp}(-\theta u_{n})/(n-1)!,$ and $p(u_{n})=\int_{0}^{\infty
}p_{\theta}(u_{n})\pi(\theta)~{{\mathrm{d}}}\theta=nu_{n}^{n-1}(u_{n}%
+1)^{-(n+1)},$ and thus\footnote{This can be obtained also from the known fact
	that in this setup the posterior distribution on $\Theta$ is a gamma
	distribution, specifically, $\Gamma(n+1,u_{n}+1),$ and so $q_{n}^{\theta_{0}%
	}=(u_{n}+1)^{n+1}\theta_{0}^{n}\operatorname{exp}(-(u_{n}+1)\theta_{0}%
	)/n!$\thinspace.}
\begin{equation}
\psi_{\theta,\theta}(n)=\pi(\theta)\int_{0}^{\infty}\frac{p_{\theta}%
	(u_{n})^{2}}{p(u_{n})}\,{{\mathrm{d}}}u_{n}=\frac{e^{-\theta}\theta^{2n}%
}{(n-1)!n!}\int_{0}^{\infty}u^{n-1}(u+1)^{n+1}\operatorname{exp}(-2\theta
u)\,{{\mathrm{{{\mathrm{d}}}}}}u. \label{eq:exppsi}%
\end{equation}

Our result is

\begin{theorem}
	\label{th:exp} Let $x_{1},x_{2},...$ be iid ${{\mathrm{Exp}}}(\theta)$
	observations, and let the prior distribution of $\theta$ be the
	${{\mathrm{Exp}}}(1)$ distribution on $\Theta=(0,\infty).$
	
	(i) For every $\theta>0$ and $n\geq1$ we have
	\[
	\psi_{\theta,\theta}(n)=\frac{\theta^{n-1/2}}{2^{n+1/2}n!\sqrt{\uppi }%
	}[(n+\theta)K_{n+1/2}(\theta)+\theta K_{n-1/2}(\theta)],
	\]
	where $K_{\nu}$ denotes the modified Bessel function of the second
	kind,\footnote{See, e.g., \cite{Abram}, Chapter 9.6.} and by \eqref{eq:expp} for every
	$\theta_{0},\theta_{1}>0$ with $(\theta_{0}+\theta_{1})/2=\theta$, we have
	\[
	\psi_{\theta_{0},\theta_{1}}(n)=\operatorname{exp}(\theta-\theta_{0}%
	)\psi_{\theta,\theta}(n)\left(  \frac{\theta_{0}\theta_{1}}{\theta^{2}%
	}\right)  ^{n}.
	\]

	(ii) For every $\theta_{0},\theta_{1}>0$ the sequence $\psi_{\theta_{0}%
		,\theta_{1}}(n)$ is strictly log-concave for $n\geq1.$
\end{theorem}

We start with two preliminary results. Put
\[
k_{n}{\;:=\;}K_{n+1/2}(\theta)
\]
and
\[
I_{n,m}:=\int_{0}^{\infty}u^{n}(u+1)^{m}e^{-2\theta u}{{\mathrm{\,{{\mathrm{d}%
}}}}}u.
\]

\begin{lemma}
	\label{l:I_n,n}For every $n\geq0$ we have
	\[
	I_{n,n}=\frac{e^{\theta}}{\sqrt{\uppi }}\frac{n!}{(2\theta)^{n+1/2}}k_{n}.
	\]
	
\end{lemma}

\begin{proof}
	The change of variable $v=2u+1$ and Equation (9.6.23) of  \cite{Abram} give
	\begin{align*}
	I_{n,n} &  =\int_{0}^{\infty}u^{n}(u+1)^{n}\operatorname{exp}(-2\theta
	u){{\mathrm{\,{{\mathrm{d}}}}}}u=\frac{e^{\theta}}{2^{2n+1}}\int_{0}^{\infty
	}((2u+1)^{2}-1)^{n}\operatorname{exp}(-\theta(2u+1))\,{{\mathrm{d}}}(2u+1)\\
	&  =\frac{e^{\theta}}{2^{2n+1}}\int_{1}^{\infty}(v^{2}-1)^{n}%
	\operatorname{exp}(-\theta v)\,{{\mathrm{d}}}v=\frac{e^{\theta}}{2^{2n+1}%
	}\frac{n!}{\sqrt{\uppi}}\left(  \frac{2}{\theta}\right)  ^{n+1/2}%
	K_{n+1/2}(\theta).
	\end{align*}
	
\end{proof}

\begin{lemma}
	\label{l:I-n-1,n+1}For every $n\geq1$ we have
	\[
	I_{n-1,n+1}=\frac{n+\theta}{n}I_{n,n}+\frac{1}{2}I_{n-1,n-1}.
	\]
	
\end{lemma}

\noindent\textit{\textbf{Proof.}}\quad%\begin{proof}
	First, the identity $u^{n-1}(u+1)^{n}=u^{n}(u+1)^{n-1}+u^{n-1}(u+1)^{n-1}$
	gives
	\begin{equation}
	I_{n-1,n}=I_{n,n-1}+I_{n-1,n-1}.\label{eq:I=I+I}%
	\end{equation}
	Second, integration by parts yields
	\begin{align*}
	2\theta I_{n,n} &  =\int_{0}^{\infty}u^{n}(u+1)^{n}2\theta\operatorname{exp}%
	(-2\theta u)\,{{\mathrm{d}}}u=\left[  u^{n}(u+1)^{n}\left(
	-\operatorname{exp}(-2\theta u)\right)  \right]  _{0}^{\infty}\\
	&  -\int_{0}^{\infty}nu^{n-1}(u+1)^{n}\operatorname{exp}(-2\theta
	u)\,{{\mathrm{d}}}u-\int_{0}^{\infty}nu^{n}(u+1)^{n-1}\operatorname{exp}%
	(-2\theta u)\,{{\mathrm{d}}}u\\
	&  =[0-0]+nI_{n-1,n}+nI_{n,n-1}%
	\end{align*}
	(we used $n>0$ for the value of the integrand at $u=0).$ By (\ref{eq:I=I+I})
	we get
	\[
	2\theta I_{n,n}=(n+n)I_{n,n-1}+nI_{n-1,n-1},
	\]
	and thus
	\begin{equation}
	I_{n,n-1}=\frac{\theta}{n}I_{n,n}-\frac{1}{2}I_{n-1,n-1}.\label{eq:n,n-1}%
	\end{equation}
	Finally, the identity $u^{n-1}(u+1)^{n+1}=u^{n}(u+1)^{n}+u^{n-1}(u+1)^{n-1}%
	+$\linebreak$u^{n}(u+1)^{n-1}$ and (\ref{eq:n,n-1}) yield
	\[
	I_{n-1,n+1}=I_{n,n}+I_{n-1,n-1}+I_{n,n-1}=\left(  1+\frac{\theta}{n}\right)
	I_{n,n}+\left(  1-\frac{1}{2}\right)  I_{n-1,n-1}.
\qquad \qquad \qed	\]
	
%\end{proof}

%We now prove Theorem \ref{th:exp}.

\begin{proof}
	[Proof of Theorem \ref{th:exp}](i) By the previous two lemmas and
	\eqref {eq:exppsi}, and setting $\psi(n):=\psi_{\theta,\theta}(n),$ we have
	\[
	\psi(n)=\frac{e^{-\theta}\theta^{2n}I_{n-1,n+1}}{(n-1)!n!}=\frac{1}%
	{\sqrt{\uppi}}\frac{\theta^{2n}}{(n-1)!n!}\left(  \frac{n+\theta}{n}\frac
	{n!}{(2\theta)^{n+1/2}}k_{n}+\frac{1}{2}\frac{(n-1)!}{(2\theta)^{n-1/2}%
	}k_{n-1}\right)  ,
	\]
	which simplifies to the claimed formula. 
	
	(ii) Let $\phi_{n}:=(n+\theta)k_{n}+\theta k_{n-1}$ and $\rho_{n}%
	:=k_{n-1}/k_{n}.$ Using the recursion
	\begin{equation}
	k_{n+1}=\frac{2n+1}{\theta}k_{n}+k_{n-1}\label{eq:recurse}%
	\end{equation}
	(by $K_{\nu+1}(u)=(2\nu/u)K_{\nu}(u)+K_{\nu-1}(u),$ which is (9.6.26) in
	\cite{Abram}) we have
	\begin{align*}
	\phi_{n+1} &  =(n+1+\theta)k_{n+1}+\theta k_{n}=(n+1+\theta)\left(
	\frac{2n+1}{\theta}k_{n}+k_{n-1}\right)  +\theta k_{n}\\
	&  =\left(  \frac{2n^{2}+3n+1}{\theta}+2n+1+\theta\right)  k_{n}%
	+(n+1+\theta)k_{n-1}\\
	&  =\left[  \left(  \frac{2n^{2}+3n+1}{\theta}+2n+1+\theta\right)
	+(n+1+\theta)\rho_{n}\right]  k_{n},\\
	\phi_{n} &  =(n+\theta)k_{n}+\theta k_{n-1}=\left[  n+\theta+\theta\rho
	_{n}\right]  k_{n},\\
	\phi_{n-1} &  =(n-1-\theta)k_{n-1}+\theta k_{n-2}=(n-1-\theta)k_{n-1}%
	+\theta\left(  k_{n}-\frac{2n-1}{\theta}k_{n-1}\right)  \\
	&  =\theta k_{n}-(n-\theta)k_{n-1}=\left[  \theta-(n-\theta)\rho_{n}\right]
	k_{n}.
	\end{align*}
	Therefore
	\begin{equation}
	R_{n}{\;:=\;}\frac{\psi(n-1)\psi(n+1)}{\psi(n)^{2}}=Q_{n}(\rho_{n}),\nonumber
	\end{equation}
	where
	\begin{equation}
	Q_{n}(\rho):=\frac{n\left[  \left(  \frac{2n^{2}+3n+1}{\theta}+2n+1+\theta
		\right)  +(n+1+\theta)\rho\right]  \left[  \theta-(n-\theta)\rho\right]
	}{(n+1)(n+\theta+\theta\rho)^{2}}.\label{eq:Q}%
	\end{equation}
	The proof that
	$Q_{n}(\rho_{n})<1\label{eq:R<1}$  
	for all $n\geq2$ and $\theta\ge 0$ is quite technical and is relegated to
	Appendix \ref{sec-a:expexp}.
\end{proof}

\textbf{General exponential prior: }If the prior on $\theta$ is a general
${{\mathrm{Exp}}}(\lambda)$ distribution for some $\lambda>0$ not necessarily
equal to $1,$ then the linear transformation $x\rightarrow x/\lambda$ reduces
it to the case of Theorem \ref{th:exp}, and so the sequence of expected
posteriors is log-concave, and thus unimodal, for $n\geq1.$

\appendix

%\section{Appendix}

The appendix is devoted to extensions of the results and to the relegated technical proofs.

\subsection{General probability models\label{sec:general}}

The results in Sections
\ref{sec:int}--\ref{sec:monot} are stated for discrete models. We now discuss
the changes needed for general models, as in Sections \ref{sec:asympt} and
\ref{sec:loggcon}.

The general setup (see, e.g., \cite{Schervish}) consists of a space of
observations $\mathcal{X}$ endowed with a measure $\nu$ and a space of
parameters $\Theta$ endowed with a measure $\mu.$ The measures $\nu$ and $\mu$
are $\sigma$-finite; for example, the counting measure in the discrete case,
and the Lebesgue measure in the continuous case. The prior probability $\Pi$
on $\Theta$ has a density function $\pi(\theta)$ with respect to $\mu,$ and
for every $\theta$ in $\Theta$ the conditional-on-$\theta$ probability
$\mathbb{P}_{\theta}$ on $\mathcal{X}$ has a density function $p_{\theta}(x)$
with respect to\footnote{These densities (which may be discrete or continuous)
	are the corresponding Radon--Nikodym derivatives.} $\nu.$ Let $\mathbb{P}%
=\int_{\Theta}\mathbb{P}_{\theta}\,\pi(\theta)\,\mathrm{d}\mu(\theta)$ be the
marginal probability on $\mathcal{X},$ with density $p(x)=$ $\int_{\Theta
}p_{\theta}(x)\pi(\theta)\,\mathrm{d}\mu(\theta)$ (for $x\in\mathcal{X}$).
Conditional on a sequence of observations $s_{n}=(x_{1},...,x_{n})$ in
$\mathcal{X}^{n}$ (endowed with the measure $\nu^{n}),$ the \emph{posterior
	density} (with respect to $\mu)$ on $\Theta$ is given at a point $\theta$ in
$\Theta$ by
$$
q_{n}^{\theta}(s_{n})=\frac{p_{\theta}(s_{n})\pi(\theta)}{p(s_{n})};
$$
the denominator is positive and finite for $\mathbb{P}$-a.e.$~s_{n}$; see
\cite{Schervish}, Theorem 1.31. All the formulas and results in the discrete
case carry over to densities in a straightforward manner; for example, the
left-hand side of (\ref{eq:martingale}) is now
\[
\mathbb{E}[q_{n+1}^{\theta_{0}}|s_{n}]=\int_{\mathcal{X}}q_{n+1}^{\theta_{0}%
}(s_{n+1})\frac{p(s_{n+1})}{p(s_{n})}\,\mathrm{d}\nu(x_{n+1}),
\]
where $s_{n+1}=(s_{n},x_{n+1}).$

The \emph{likelihood ratio order} given in Section \ref{sec:sing} is defined
for general random variables $x$ and $y$ as follows: $y\geq_{{\mathrm{lr}}}x$
if
\begin{equation}
{\mathbb{P}}(x\in A)\,{\mathbb{P}}(y\in B)\geq{\mathbb{P}}(x\in
B)\,{\mathbb{P}}(y\in A)\label{eq:SSLR}%
\end{equation}
for any two sets $A$ and $B$ in $\mathbb{R}$ such that $A\leq B,$ which means
that $a\leq b$ for all $a\in A$ and $b\in B$; see \cite{Shak}, (1.C.3).
Everything that is stated about this order relation in Section \ref{sec:sing},
including its relation to the stochastic order (whose definition in terms of
expectations of increasing functions $f$ remains the same) continues to hold.

Part (i) of Proposition \ref{prop:p1p2} is now stated as follows. Let
$P_{1}\ll P_{2}$ be two probability measures on a space ${{\mathcal{S}}}$ with
corresponding densities $p_{1}$ and $p_{2}$ with respect to the underlying
measure $\nu$ on $\mathcal{S}.$ Put $r(s):=p_{1}(s)/p_{2}(s)$ for every $s$
(again, when the ratio is $0/0,$ which has probability $0$ under both $P_{1}$
and $P_{2},$ define $r(s)$ arbitrarily); then
\[
{\mathcal{L}}_{P_{1}}(r)\geq_{{{\mathrm{lr}}}}{\mathcal{L}}_{P_{2}}(r),
\]
i.e., $P_{1}\circ r^{-1}\geq_{{{\mathrm{lr}}}}P_{2}\circ r^{-1}.$

\begin{proof}
	Let $A,B\subset\lbrack0,\infty)$ be such that $A\leq B;$ i.e., there is $c$
	such that $\sup A\leq c\leq\operatorname{inf}B.$ We need to show that
	\begin{equation}
	P_{1}(r^{-1}(B))P_{2}(r^{-1}(A))\geq P_{1}(r^{-1}(A))P_{2}(r^{-1}%
	(B)).\label{eq:lr-AB}%
	\end{equation}

	If $c=0$ then $A=\{0\},$ and so for every $s\in r^{-1}(A)$ we have $r(s)=0,$
	and thus $p_{1}(s)=0;$ therefore $P_{1}(r^{-1}(A))=\int_{r^{-1}(A)}%
	p_{1}(s)\,{\mathrm{d}}\nu(s)=0,$ and (\ref{eq:lr-AB}) holds since its
	right-hand side is $0.$
	
	Let thus $c>0.$ For every $s\in r^{-1}(A)$ we have $r(s)\leq c,$ and thus
	$p_{1}(s)\leq cp_{2}(s);$ integrating over $r^{-1}(A)$ gives
	\[
	P_{1}(r^{-1}(A))\leq cP_{2}(r^{-1}(A)).
	\]
	Similarly, for every $s\in r^{-1}(B)$ we have $r(s)\geq c,$ and thus
	$p_{2}(s)\leq(1/c)p_{1}(s),$ which, integrating over $r^{-1}(B),$ gives
	\[
	P_{2}(r^{-1}(B))\leq\frac{1}{c}P_{1}(r^{-1}(B)).
	\]
	Multiplying the two inequalities yields (\ref{eq:lr-AB}).
\end{proof}

Part (ii) of Proposition \ref{prop:p1p2}, i.e.,%
\[
{{\mathcal{L}}}_{\theta}(q^{\theta})\geq_{{{\mathrm{lr}}}}{{\mathcal{L}}%
}(q^{\theta}),
\]
follows now from part (i) when $\mathbb{P}_{\theta}\ll\mathbb{P}$ (see
\cite{Bill}, (34.15)), which clearly holds in our setup where $\Theta
\subseteq\mathbb{R}$ is an interval, the prior $\pi$ is positive on $\Theta,$
and the densities $p_{\theta}$ are continuous in $\theta.$

Finally, we note that while posterior probabilites are always bounded from
above by $1$ when the space of parameters $\Theta$ is discrete, posterior
densities need \emph{not} be bounded in general (as seen in Sections
\ref{sec:asympt}--\ref{sec:loggcon} and Appendix \ref{sec-a:asymp} below).

\subsection{Asymptotic analysis\label{sec-a:asymp}}

Section \ref{sec:asympt} deals with one-dimensional exponential families of
distributions. We expect the analysis to extend to more general setups, in
particular, to multidimensional exponential families. Such a family is given
by densities $p_{\theta}(x)=\exp(\eta(\theta)\cdot T(x)-A(\eta(\theta))-B(x))$
(for $\theta\in\Theta$ and $x\in\mathcal{X)},$ where $d\geq1,$ $\eta
:\Theta\rightarrow\mathbb{R}^{d},$ $T:\mathcal{X\rightarrow\mathbb{R}}^{d},$
$A:\mathbb{R}^{d}\rightarrow\mathbb{R},$ and $B:\mathcal{X}\rightarrow
\mathbb{R}$ (in the \textquotedblleft canonical" representation, $\eta$ is the
identity and $\Theta\subseteq\mathbb{R}^{d}).$

\subsubsection{{Generalizing Proposition \ref{p:theta-bar}}}

Consider first the useful reduction to the $\theta_{1}=\theta_{0}$ case. Given
two densities $p_{0}$ and $p_{1}$ on the space $\mathcal{X}$ (with respect to
a positive $\sigma$-finite measure $\nu),$ define their \emph{normalized
	geometric average} (\emph{NGA) }to be the density $r$ on $\mathcal{X}$ given
by
\[
r(x)%
%TCIMACRO{\TeXButton{:=}{{\;:=\;}}}%
%BeginExpansion
{\;:=\;}%
%EndExpansion
\frac{p_{0}(x)\,p_{1}(x)}{\int_{{{\mathcal{X}}}}\sqrt{p_{0}(x)\,p_{1}%
		(x)}\,{{\mathrm{d}}}\nu(x)}%
\]
for every $x\in{{\mathcal{X}}}$. Using (\ref{eq:psi}), Proposition
\ref{p:theta-bar} readily generalizes to

\begin{proposition}
	\label{p:theta2}Let $(p_{\theta})_{\theta\in\Theta}$ be a family of densities
	and let $\pi$ be a prior density on $\Theta$ with $\pi(\theta)>0$ for all
	$\theta\in\Theta.$ Let $\theta_{0},\theta_{1}\in\Theta$ and let $r$ be the
	normalized geometric average of $p_{\theta_{0}}$ and $p_{\theta_{1}}.$ If
	$r=p_{\theta_{2}}$ for some $\theta_{2}\in\Theta$ then%
	\begin{equation}
	\psi_{\theta_{0},\theta_{1}}(n)=\frac{\pi(\theta_{0})}{\pi(\theta_{2})}%
	\psi_{\theta_{2},\theta_{2}}(n)\,w^{n},\label{eq:reduction}%
	\end{equation}
	where $w$ is given by (\ref{eq:w}).
\end{proposition}

For multidimensional exponential families, closure under normalized geometric
averages (\emph{NGA-closure} for short) is easily seen to amount to the
convexity of the set of  \textquotedblleft natural parameters" $\eta
(\Theta)$ (because $\eta(\theta_{2})$ must equal $(\eta(\theta_{0}%
)+\eta(\theta_{1}))/2$). In the one-dimensional case with $\Theta
\subseteq\mathbb{R}$ and $\eta^{\prime}>0$ the convexity of $\eta(\Theta)$
follows from the convexity of $\Theta.$

For families of distributions that are not NGA-closed, one may consider for
each $p_{\theta_{0}}$ and $p_{\theta_{1}}$ that member of the family,
$p_{\theta_{2}}$ with $\theta_{2}\in\Theta,$ that is closest to their
normalized geometric average $r$ in terms of the Kullback--Leibler distance.
Indeed, this $p_{\theta_{2}}$ yields the \textquotedblleft leading exponential
term" in (\ref{eq:psi}) as $n\rightarrow\infty$ (cf. the Laplace method, which
is used in the Bernstein--von Mises result as well; this explains why
$p_{\theta_{2}}$ must belong to the family and $\pi(\theta_{2})>0$). The
simplicity of the relation (\ref{eq:reduction}) is however lost when
$p_{\theta_{2}}\neq r$.

Interestingly, for single-parameter families of densities ${\mathcal{F=\{}%
}p_{\theta}:\theta\in\lbrack a,b]\}$ (for $a<b)$ where the dependence on
$\theta$ is continuous and there is identifiability ($p_{a}\neq p_{b}$
suffices), $\mathcal{F}$ is NGA-closed if and only if $\mathcal{F}$ is a
one-dimensional exponential family (\ref{eq:exp-fam}) on $[a,b].$ Indeed,
NGA-closure implies by continuity closure with respect to all normalized
geometric weighted averages; i.e., for every $p_{\theta_{0}},p_{\theta_{1}}%
\in\mathcal{F}$ and every $\lambda\in\lbrack0,1]$ there is $\tau(\lambda
)\in\lbrack a,b]$ such that $c(\lambda)(p_{\theta_{0}})^{1-\lambda}%
(p_{\theta_{1}})^{\lambda}=p_{\tau(\lambda)},$ for the appropriate
normalization constant $c(\lambda).$ Take $\theta_{0}=a$ and $\theta_{1}=b;$
then $\tau(0)=a$ and $\tau(1)=b;$ the function $\tau$ is easily seen to be
continuous and one-to-one (because $p_{a}\neq p_{b}$ and so all the
$p_{\tau(\lambda)}$ are distinct), and thus $\tau$ is strictly increasing and
its range is the whole interval $[a,b].$ From $\log p_{\tau(\lambda)}=\log
c(\lambda)+(1-\lambda)\log p_{a}(x)+\lambda\log p_{b}(x)$ we thus get%
\[
\log p_{\theta}(x)=\tau^{-1}(\theta)(\log p_{b}(x)-\log p_{a}(x))+\log
p_{a}(x)+\log c(\tau^{-1}(\theta))
\]
for every $\theta\in\lbrack a,b]$ and $x\in\mathcal{X},$ which yields the
one-dimensional family (\ref{eq:exp-fam}) with $\eta=\tau^{-1},$ $T(x)=\log
p_{b}(x)-\log p_{a}(x),$ $A(\eta)=-\log(c(\eta)),$ and $B(x)=-\log p_{a}(x)$.

\subsubsection{{Generalizing Theorem \ref{th:asymp}}}

Consider the asymptotic result of Theorem \ref{th:asymp} (i) when $\theta
_{1}=\theta_{0},$ but now for a multivariate $d$-dimensional exponential
family as above. We expect the rate $\sqrt{n}$ to become $n^{d/2}$ (with
appropriate constants). Indeed, in Step 1 the change of variable
$\vartheta=(nI(\widehat{\theta}_{n}))^{1/2}(\theta-\widehat{\theta}_{n})$ in
the Bernstein--von Mises theorem (where $nI$ is now a $d\times d$ matrix)
involves a Jacobian, and so yields a factor of $n^{d/2};$ in Step 4, the bound
on $q_{n}^{\theta_{0}}$ becomes $Cn^{d/2},$ as it involves a $d$-dimensional
Gaussian integral (cf. the $d$-dimensional Laplace method).

 We conjecture that the analysis extends beyond exponential families (for
 instance, in setups where the Bernstein--von Mises result
 applies).
 
\subsubsection{Proof of Theorem \ref{th:asymp}: details for Step
	1\label{sec-a:step1}}

We provide here technical details for Step 1 of the proof of Theorem
\ref{th:asymp}\textbf{.}

First, we show that all $7$ regularity conditions of Theorem 7.89 in
\cite{Schervish} (the Bernstein--von Mises result that we use) hold for an
exponential family (\ref{eq:exp-fam-canon}). Indeed, this is immediate for
conditions 1--4 (see the discussion in Section \ref{sec:asympt}). Condition 5:
in our setup $\lambda_{n}=1/(nA^{\prime\prime}(\hat{\theta}_{n}))$, and so
$\lambda_{n}\overset{{{\mathbb{P}}}_{0}}{\longrightarrow}0$ follows from
$A^{\prime\prime}(\hat{\theta}_{n})\overset{{{\mathbb{P}}}_{0}%
}{\longrightarrow}I(\theta_{0})>0.$ Condition 6 is equivalent to the
existence, for each $\delta>0$, of a positive $K(\delta)$ such that
\[
\lim_{n\rightarrow\infty}{{\mathbb{P}}}_{0}\left(  \sup_{\theta\notin%
	\lbrack\theta_{0}-\delta,\,\,\theta_{0}+\delta]}\,\frac{1}{nA^{\prime\prime
	}(\hat{\theta}_{n})}\sum_{i=1}^{n}[\log p_{\theta}(x_{i})-\log p_{\theta_{0}%
}(x_{i})]<-K(\delta)\right)  =1.
\]
To see this note that $(1/n)\sum_{i=1}^{n}[\log p_{\theta}(x_{i})-\log
p_{\theta_{0}}(x_{i})]=A(\theta_{0})-A(\theta)+(\theta-\theta_{0})\overline
{t}_{n}$, and since $\overline{t}_{n}\overset{{{\mathbb{P}}}_{0}%
}{\longrightarrow}A^{\prime}(\theta_{0})$, it suffices to show that
$A(\theta_{0})-A(\theta)+(\theta-\theta_{0})A^{\prime}(\theta_{0})<0$. The
latter inequality follows readily from the strict convexity of $A$. Finally,
condition 7 requires that for each $\varepsilon>0$ there exist $\delta>0$
such that
\begin{multline*}
\lim_{n\rightarrow\infty}{{\mathbb{P}}}_{0}\left(  \sup_{\theta\in
	\lbrack\theta_{0}-\delta,\theta_{0}+\delta]}\left\vert 1+{\lambda_{n}}%
\sum_{i=1}^{n}\frac{\mathrm{d}^{2}}{\mathrm{d}\theta^{2}}\log p_{\theta}%
(x_{i})\right\vert <\varepsilon\right)  \\
=\lim_{n\rightarrow\infty}{{\mathbb{P}}}_{0}\left(  \sup_{\theta\in
	\lbrack\theta_{0}-\delta,\theta_{0}+\delta]}\left\vert 1-{n\lambda_{n}%
}A^{\prime\prime}(\theta)\right\vert <\varepsilon\right)  =1.
\end{multline*}
We have $n\lambda_{n}\overset{{{\mathbb{P}}}_{0}}{\longrightarrow}%
1/I(\theta_{0})=1/A^{\prime\prime}(\theta_{0})$ implying, together with the
continuity of $A^{\prime\prime}$, that there exists $\delta$ such that
${n\lambda_{n}}A^{\prime\prime}(\theta)=A^{\prime\prime}(\theta)/A^{\prime
	\prime}(\hat{\theta}_{n})\in\lbrack1-\varepsilon,1+\varepsilon]$ for
$\theta\in\lbrack\theta_{0}-\delta,\theta_{0}+\delta]$ holds with probability
converging to $1$ as $n\rightarrow\infty$.

Second, to translate the density $f_{\vartheta|s_{n}}$ of $\vartheta
=\sqrt{nI(\hat{\theta}_{n})}(\theta-\hat{\theta}_{n})$ given $s_{n}$ to that
of $\theta$ given $s_{n}$---which is the posterior $q_{n}^{\theta}(s_{n}%
)$---we need to divide the latter by $\sqrt{nI(\hat{\theta}_{n})}$.

\subsection{The normal case\label{sec-a:normnorm}}

We provide here technical details for the proof of Theorem \ref{th:normal} (i)
in Section \ref{sec:normnorm}.

\begin{proof}
	[Proof of Theorem \ref{th:normal} (i)]Recall that the distribution $p_{\theta
	}(u_{n})$ of the sufficient statistic $u_{n}=\sum_{i=1}^{n}x_{i}$ given
	$\theta$ is ${{\mathcal{N}}}(n\theta,n\sigma^{2}),$ and its marginal
	distribution $p(u_{n})$ is ${{\mathcal{N}}}(0,n^{2}+n\sigma^{2}).$ Therefore
	$\psi(n)\equiv\psi_{\theta,\theta}(n)=\int p_{\theta}^{2}(u_{n})/p(u_{n}%
	)\,{{\mathrm{d}}}u_{n}=a_{n}I_{n},$ where $I_{n}=\int_{-\infty}^{\infty
	}\operatorname{exp}\left(  -h_{n}(u)\right)  {{\mathrm{d}}}u$ for
	\[
	h_{n}(u)=2\frac{(u-n\theta)^{2}}{2n\sigma^{2}}-\frac{u^{2}}{2(n^{2}%
		+n\sigma^{2})},
	\]
	and $a_{n}$ is the constant
	\[
	a_{n}=\left(  \frac{1}{\sqrt{2\uppi}\sqrt{n}\,\sigma}\right)  ^{2}\left(
	\frac{1}{\sqrt{2\uppi}\sqrt{n^{2}+n\sigma^{2}}}\right)  ^{-1}=\frac
	{\sqrt{n^{2}+n\sigma^{2}}}{\sqrt{2\uppi}\,n\sigma^{2}}.
	\]
	The function $h_{n}(u)$ is a quadratic in $u,$ namely, $h_{n}(u)=\left(
	b_{n}u-c_{n}\right)  ^{2}-d_{n}$ for
	\[
	b_{n}=\sqrt{\frac{2n+\sigma^{2}}{2n\sigma^{2}(n+\sigma^{2})}}%
	\,{{\text{\ \ and\ \ \ }}}d_{n}=\frac{\theta^{2}n}{2n+{\sigma}^{2}}%
	\]
	(the value of $c_{n}$ will not matter), which yields $I_{n}=(\sqrt
	{\uppi}/b_{n})\operatorname{exp}(d_{n}).$ Substituting in $\psi(n)=a_{n}I_{n}$
	gives (\ref{eq:normal}).
\end{proof}

\subsection{The exponential case\label{sec-a:expexp}}

We prove here the final inequality in the proof of Theorem
\ref{th:exp} (ii) in Section \ref{sec:expexp}, namely $Q_{n}(\rho_{n})<1$ for
all $n\geq2,$ and $\theta \ge 0$; the function $Q_{n}(\rho)$ is defined in (\ref{eq:Q}) and
$\rho_{n}=k_{n-1}/k_{n}\equiv K_{n-1/2}(\theta)/K_{n+1/2}(\theta).$

\begin{proof}
	 	We will use the following bounds on $\rho_{n}$:
	\begin{equation}
	\eta_{n}:=\frac{\theta}{n+\frac{1}{2}+\sqrt{\left(  n-\frac{3}{2}\right)
			^{2}+\theta^{2}}}<\rho_{n}\leq\theta;\label{eq:L}%
	\end{equation}
	see \cite{Segura}: the lower bound is from (34), and the upper bound from
	(33). Taking the derivative of $Q_{n}$ yields
	\begin{multline*}
	Q_{n}^{\prime}(\rho)\\
	=-{\frac{2{n}^{5}+\left(  2\theta+3\right)  {n}^{4}+\left(  2\theta^{2}%
			+\theta+1+2{\theta}^{2}\rho-\theta\rho\right)  {n}^{3}+\left(  2\theta
			^{2}+2\theta^{2}\rho-\theta\rho\right)  n^{2}+\,\left(  \theta^{2}+\theta
			^{2}\rho\right)  n}{(n+1)\left(  n+\theta+\theta\rho\,\right)  ^{3}\theta};}%
	\end{multline*}
	in the range $0\leq\rho\leq\theta$ the coefficients of all powers of $n$ in
	the numerator are positive (use $\theta\rho\leq2\theta^{2}),$ and so $Q_{n}$
	is strictly decreasing in $\rho$ there. Since $0\leq\eta_{n}<\rho_{n}%
	\leq\theta$ by (\ref{eq:L}), it follows that
$
	R_{n}=Q_{n}(\rho_{n})<Q_{n}(\eta_{n}),
$
	and so it suffices to show that $Q_{n}(\eta_{n})<1$ for all $n\geq2.$
	Computing $Q_{n}(\eta_{n})$ by substituting (\ref{eq:L}) in (\ref{eq:Q})
	yields%
	\[
	1-Q_{n}(\eta_{n})=\frac{A+B\sqrt{C}}{D},
	\]
	where we put $n=m+2$ (and so $m\geq0$ when $n\geq2),$
	\begin{align*}
	A &  =8\,{m}^{4}+72\,{m}^{3}+\left(  4\,{\theta}^{2}-16\,\theta+230\right)
	{m}^{2}+\left(  8\,{\theta}^{3}+4\,{\theta}^{2}-56\,\theta+302\right)  m\\
	&  +(8\,{\theta}^{4}+20\,{\theta}^{3}+2\,{\theta}^{2}-48\,\theta+132),\\
	B &  =4\,{m}^{3}+30\,{m}^{2}+\left(  -4\,{\theta}^{2}+74\right)
	m+(-4\,{\theta}^{3}-10\,{\theta}^{2}+60),\,\,
	C  =4m^{2}+4m+(4\theta^{2}+1),
	\end{align*}
	and the denominator $D$ is positive. We claim that in the range $m\geq0$ and
	$\theta\geq0$ we have $A>0$ and $E:=A^{2}-B^{2}C>0,$ and thus $A+B\sqrt{C}>0$
	(immediate when $B\geq0;$ when $B<0,$ use $A+B\sqrt{C}=E/(A-B\sqrt{C})$).
	Indeed, the coefficients of the powers of $m$ in $A$ and $E$ are positive for all $\theta \ge 0$, as shown by direct calculations that we omit.
Thus $1-Q_{n}(\eta_{n})>0$ for all $n\geq2$ and $\theta\geq0$. 
\end{proof}

\noindent\textbf{Acknowledgment} We thank Marco Scarsini for many useful discussions. We are grateful to the referees, the associate editor, and the editor for many constructive comments, and in particular for encouraging us to consider exponential families and the Bernstein-von Mises Theorem.

\end{document}